\documentclass[11pt,reqno]{amsart}

\usepackage{amscd,amsfonts,amssymb}
\usepackage{graphicx}
\usepackage[english]{babel}
\usepackage{epsfig}

\copyrightinfo{2000}{American Mathematical Society}

\numberwithin{equation}{section}

\theoremstyle{definition}

\theoremstyle{remark}

\theoremstyle{cor 3.} \theoremstyle{cor 2.}
\theoremstyle{Proposition 2.} \theoremstyle{Proposition 3.}
\newtheorem{cor 2.}{Corollary 2.}
\newtheorem{cor 3.}{Corollary 3.}
\newtheorem{utv 2.}{Proposition 2.}
\newtheorem{utv 3.}{Proposition 3.}

\author{Olga ~S.~Rozanova}
\address{Mechanics and Mathematics Faculty, Moscow State University, Leninskie Gory 1, Moscow, Russia.} \email{ rozanova@mech.math.msu.su}

\title[Smooth solutions to the multidimensional
balance laws ]
 {Classes of smooth solutions to the multidimensional
balance laws of gas dynamic type on the riemannian manifolds}

\thanks{}


\date{July 11, 2002 }

\begin{document}

\begin{abstract}
The paper is devoted to the special classes of solutions to
multidimensional balance laws of gas dynamic type. In the velocity
field for the solutions of such class the time and space variables
are separated. The simplest case is the solution with the linear
profile of velocity in the Euclidean space.
%Also we apply the theory of solutions with linear profile of velocity
%to the problem of tracing of the typhoon eye.
 \end{abstract}

\maketitle

%\begin{document}

%%% ----------------------------------------------------------------------

%\documentclass[reqno,12pt]{article}
%\usepackage[english]{babel}
%\usepackage{graphicx,graphics}
%\usepackage{textcomp}
%\usepackage{amsmath}
%\usepackage{amssymb}

%\usepackage [dvips]{graphics}
%\usepackage [dvips]{graphicx}
%\usepackage[pdftex]{graphicx}

%\textheight= 227mm
%\textwidth=160mm
%\hoffset=-1.7cm

\parindent=0,5cm

\section {Balance laws on the riemannian manifold}
\subsection{Preliminaries}

Let $\Sigma$ be a smooth $n$ -- dimensional riemannian manifold with
local curvilinear orthogonal coordinates ${x}^i,\,i=1,..,n.$
Consider the following system of nonlinear partial differential
equations
$$(R^i({\bf U},{\bf x},t))'_t +
\nabla _j q^{ij}
({\bf U},{\bf x},t) + Q^i ({\bf U},{\bf x},t) = 0,\quad
i,j=1,\dots,n.\eqno(1.1)$$

Here
${\bf U}=(u^1,\dots,u^n),$ ${\bf R}=(R^1,\dots,R^n)$
and
${\bf Q}=(Q^1,\dots,Q^n)$
are the vector-functions on the tangent bundle of
$\Sigma$ ($u$ is unknown),
$q^{ij}$ is a smooth tensor field, $\nabla _j
q^{ij}$ is a divergency of tensor. All differential operations
with respect to space variables are performed in this coordinate
system by means of smooth metric tensor
$ g_{ij} $ of the manifold
$\Sigma$, the variable $ t$ denotes the time, as usual.
Recall that the covariant derivative of any contravariant
vector ${\bf
X}(X^1,\dots,X^n)$ is the tensor $\nabla_j
X^i=\frac{\partial X^i}{\partial x^j}+\Gamma_{jk}^i X^k,\,
i,j,k=1,\dots,n ,$ where $\Gamma_{ij}^k\, $ are the Christoffel
symbols.

If the tensor field $q^{ij}$ is sufficiently smooth, then the system
(1.1) can be re-written in the form
$$G^i_j({\bf U},{\bf x},t)(u^i)'_t +
 F^{ik}_j({\bf U},{\bf x},t)\nabla _k
u^j + Q^i ({\bf U},{\bf x},t) = 0,\quad
i=1,\dots,n,\eqno(1.2)$$
with new tensor fields $G^i_j$ and $F^{ik}_j,$
or
$$
A^0({\bf U}, {\bf x}, t )\frac{\partial {\bf U}}{\partial
t}+\sum\limits^n_{k=1} A^k({\bf U}, {\bf x}, t )
\frac {\partial {\bf U}}{\partial {x}^k}+\bar {\bf Q}
({\bf U},{\bf x},t)=0,
\eqno (1.3)
$$
where $A^0=G^i_j,\,A^k=F^{ik}_j,\,\bar Q^i=Q^i+F^{ik}_m
\Gamma^m_{kj}u^j.$

Remark that the system (1.2) includes only covariant
derivatives, which coincide with the usual partial ones only
in the case of the euclidean space.  Nevertheless, the
representation of the covariant derivatives through the partial
ones and the Christoffel symbols \cite{NovikovFomenko}
adds  the derivatives free terms only.

We restrict ourself by the systems of gas dynamic type,
namely, by the system of the Euler equation with the right
hand sides of special form.

The important fact is the possibility to reduce these balance laws
to the symmetric hyperbolic form.

Recall that  the system of form
(1.3) for the unknown
function ${\bf U}(t,{\bf x})$ is called
symmetric hyperbolic if the  matrices $A^j(
{\bf U},{\bf x},t), j=0,...,n,$ are symmetric and the matrix
$A^0({\bf U},{\bf x},t) $ is positively
defined in addition.
 The following classical result is
well known for the euclidean space ${\mathbb R}^n$ (\cite {Kato},
\cite{VolpertKhudiaev}, particular cases in \cite{GardingLeray},
\cite{Lax}, \cite{Majda}, \cite{Crandall}):  if the matrices $A^j(
{\bf p},{\bf x},t)$ and the right hand side $\bar{\bf Q}( {\bf
p},{\bf x},t)$ smoothly depend on their arguments, have the
continuous bounded derivatives up to order $m+1$ with respect to the
variables $( {\bf p},{\bf x})$ under bounded $\bf p$, and the
initial data ${\bf U}(0,{\bf x})$ and the function $\bar{\bf Q}({\bf
0},{\bf x},t)$ belong to the Sobolev class $H^m({\mathcal
R}^n),\,m>1+n/2,$ then the corresponding Cauchy problem has locally
in time unique solution in the class $\cap _{j=0}^m
C^j([0,T);H^{m-j}({\mathcal R}^n))$.  Moreover, $$\lim\limits_{t\to
T-0}\sup(\|{\bf U}\|_{L^\infty}+ \|\nabla_{\bf x} {\bf
U}\|_{L^\infty})=+\infty.$$ (Remark that on the initial data less
strict requirements may be imposed, for example, the requirement of
belonging uniformly local in space to $H^m({\mathcal R}^n).$
\cite{Chemin}, \cite{VolpertKhudiaev}).

As usual in the problem arising from physics the coefficients
of system (1.3) depend only on the solution $\bf U.$ The fact
simplifies significantly the formulation of the result.  Namely, for
the local in time existence of the Cauchy problem in the described
class it is sufficiently
to require the implementation of the condition $\bar{\bf Q}({\bf
0})={\bf 0}$ besides of the smoothness of coefficients.

If the manifold $\Sigma$ is diffeomorfic to   ${\mathcal R}^n,$ then
the theorem on the existence of smooth local in time solution can be
transmitted to the case.

The author don't know the analogous result for the
case of an arbitrary riemannian manifold.

\subsection{The system of the Euler equations}

Further we shall consider the following system given on the
riemannian manifold $\Sigma$ in $n$ dimensions with the metric tensor
$g_{ij}$ (below we denote it (E1--E3))
$$ \rho (\partial_t {\bf V}+({\bf V},{\bf \nabla})\,{\bf V})= -{\bf\nabla} p + \rho {\bf
F}(\rho, p, {\bf V}, {\bf x}, t),\eqno(1.2.1) (E1)$$
$$\partial_t \rho + div\, (\rho{\bf
V})=0\eqno(1.2.2) (E2)$$
$$ \partial_t p +({\bf V},{\bf\nabla}
p)+\gamma p\, div\,{\bf V}=0.\eqno(1.2.3) (E3)$$
Here $ \rho,\, p, \,
{\bf V}$ are the density, the pressure and the velocity vector,
respectively.  The constant $\gamma$ is the heat ratio, $\gamma>1,$
moreover, in the physical case $\gamma\le 1+\frac{2}{n}.$

The right-hand side $\rho {\bf F}$ is the forcing term, it
is supposed to be  smooth with respect to all arguments.

For example, it may describe the friction if
${\bf F}=-\mu {\bf V}
|{\bf V}|^\sigma$ with the nonnegative function $\mu(t,{\bf x})$ and
the constant $\sigma\ge 0$,
or the Coriolis term. In the last case
${\bf F}=[{\bf V\times \omega}], \,{\bf\omega}$  is a constant
vector, corresponding to the vector of angular velocity (for $n=3$)
or
 ${\bf F(V)}=l{\bf V}_\bot,\, l=l({\bf x},t)$ is the Coriolis
parameter, $V_\bot^i=e_{ij}V^j,\,e_{ij}$ is the skew-symmetric
discriminant tensor of the surfaces (for $n=2$).

\section{The Euclidean space. Solutions with the linear profile of
velocity and their properties}

We consider the classical solutions to
(E1--E3) with the density sufficiently quickly decreasing
as $|{\bf x}|\to \infty$ to guarantee a convergency of the
integral $\int \limits_{{\mathbb R}^n}\rho |{\bf x}|^2 d{\bf x}$ (so
called solutions with a finite moment).

As well known, solutions to the Cauchy problem for system
(E1--E3) may lose the initial smoothness for a finite time.
Sometimes there is a possibility to estimate the time of a
singularity formation from above (see, f.e., \cite{RozDU} and
references therein). Moreover, in the case
${\bf F}=0$ the singularity appears in any solution with
compactly supported initial data (f.e., \cite
{Xin}).

At the same time, there exist some nontrivial
classes of globally smooth solutions. In the section our special
interest will be {\it the solutions with linear profile of velocity}
$${\bf V}=A(t){\bf r}+{\bf b}(t),\eqno(2.1)$$ where $A(t)$ and ${\bf
b}(t)$ are a matrix $(n\times n)$ and an $n$ - vector, dependent on
time, ${\bf r}$ is a radius-vector of the point.

\subsection{Construction of solutions with linear profile of
velocity and integral functionals}
\subsubsection{The system on the plane ($n=2$)}

Now we construct solutions with linear profile of
velocity for the important case
${\bf F}= L{\bf V},$ with matrix $L=\left(\begin{array}{cc}
-\mu&-l\\l&-\mu\end{array}\right),\, $
$ \mu=const\ge 0,\, l=const.$
This right-hand side describes in a simplest way the Coriolis
force and the Rayleigh friction in the meteorological model
neglecting the vertical processes.

We consider only smooth solutions to (E1--E3)
with the density (and pressure) vanishing at $|{\bf x}|\to\infty$
rather quickly to guarantee the convergency of all integrals
involved (whereas the velocity components may even
grow). For the solution the total mass
$${\mathcal M}=\int\limits_{{\mathbb
 R}^2} \rho d{\bf x}$$ is conserved,
the total energy $${{\mathcal E}}(t)=\int\limits_{{\mathbb
 R}^2} \left( \frac{\rho |{\bf
V}|^2)}{2}+\frac{p}{\gamma-1} \right) \,d{\bf x}={\mathcal
E}_k(t)+{\mathcal E}_p(t)$$ is nonincreasing function. Let us
involve integral
 functionals $$G(t)=\frac{1}{2}\int \limits_{{\mathbb R}^2}\rho|{\bf
r}|^2\,d{\bf x},$$
$$
N_i(t)=\int\limits_{{\mathbb
 R}^2} X_{1i}\rho\,d{\bf x},\,
i=1,2,\,$$
$$
I_i(t)=\int\limits_{{\mathbb
 R}^2} V_i \rho\,d{\bf x},\,
i=1,2,\,$$
$$
F_i(t)=\int\limits_{{\mathbb
 R}^2}({\bf V, X}_i)\rho\,d{\bf x},\,
i=1,2,\,$$ where ${\bf X}_1={\bf r}=(x,y),\,{\bf
X}_2={\bf r}_\bot=(y,-x).$
Note that $G(t)>0.$

For the {\it classical} (and for the piecewisely smooth)
solutions to
(E1 -- E3), the following relations hold \cite{Roz98FAO} (see also
\cite{Gordin}):
$$ G'(t)=F_1(t),\eqno(2.1.1)$$
$$N_i'(t)=I_i(t),\,i=1,2,\eqno(2.1.2)$$
$$I_1'(t)=-l I_2(t)-\mu I_2(t),\eqno(2.1.3)$$
$$I_2'(t)=l I_1(t)-\mu I_2(t),\eqno(2.1.4)$$
$$ F_1'(t)=2(\gamma-1){\mathcal E}_p(t)+2 {\mathcal E}_k(t)
- l F_2(t)-\mu F_1(t),\eqno(2.1.5) $$
$$F_2'(t)=l F_1(t)-\mu F_2(t),\eqno(2.1.6)$$
$${\mathcal
E}'(t)= -2\mu {\mathcal E}_k(t).\eqno(2.1.7) $$

The result can be obtained by means of the Green' formula.

Generally speaking, the functions ${\mathcal E}_p(t)$ and ${\mathcal
E}_k(t)$ cannot be expressed through $G(t)$, $F_1(t)$, $F_2(t)$,
$I_1(t)$, $I_2(t)$, $N_1(t)$, $N_2(t).$  But if we choose the
velocity field with the linear profile $ {\bf V}= A(t){\bf r}+{\bf
b}(t),$ we obtain the {\it closed} system of ODE for finding the
coefficients of the matrix $A(t)$ and of the components of vector
${\bf b}(t)$.  In general case it is rather complicated  and can be
solved only numerically.

More simple result can be obtained if
$$
{\bf V}= \alpha(t){\bf r}+\beta(t){\bf r_\bot}+
{\bf b}(t)=$$
$$=\left(\alpha(t)I+\beta(t)
\left(\begin{array}{cc}
0&1\\-1&0\end{array}\right)\right){\bf r}+
\left(\begin{array}{c}
b_1(t)\\b_2(t)\end{array}\right),\eqno(2.1.8)
$$
with the identity matrix $I.$
It is easy to see that in this case
$$
I_1(t)=\alpha(t) N_1(t)+\beta(t)N_2(t)+b_1(t){\mathcal M}, $$
$$
I_2(t)=\alpha(t) N_2(t)-\beta(t)N_1(t)+b_2(t){\mathcal M}, $$
$$
F_1(t)= 2\alpha(t) G(t)+b_1(t)N_1(t)+b_2(t)N_2(t), $$
$$  F_2(t)= 2\beta(t) G(t)+b_2(t)N_1(t)-b_1(t)N_2(t),$$
$$
{\mathcal E}_k(t)=(\alpha^2(t)+\beta^2(t)) G(t)
+\alpha(t)b_1(t)N_1(t)+\beta(t)b_1(t)N_2(t)+$$$$
+\alpha(t)b_2(t)N_2(t)+\beta(t)b_2(t)N_1(t)+
\frac{1}{2}(b_1^2(t)+b_2^2(t)){\mathcal M},
$$
$${\mathcal E}_p(t)= K
G^{1-\gamma}(t), K={\mathcal E}_p(0) G^{\gamma-1}(0).$$

 In that way, all functions involved in
the system (2.1.1--2.1.7) are expressed through $G(t),$
$N_1(t),$ $N_2(t),$
$\alpha(t),$ $ \beta(t),$ $ b_1(t),$  $ b_2(t).$

\subsubsection{$A(t)$ of special form, ${\bf b}(t)=0$}

The case is simplest. Namely, $F_1(t),$ $F_2(t),$ $E_k(t),$ $E_p(t)$
can be expressed through $G(t),\,\alpha(t),\,\beta(t).$ Therefore,
if we denote for the convenience $G_1(t)= 1/G(t)$, then we obtain
the system
$$ G_1'(t)=-2\alpha(t)
G_1(t), $$ $$ \beta'(t)=\alpha(t)(l-2\beta (t)) - \mu
\beta(t),\eqno(2.1.9)$$ $$
\alpha'(t)=-\alpha^2(t)+\beta^2(t)-l\beta(t)-\mu\alpha(t)+
(\gamma-1)K G_1^{\gamma}(t).
$$

For $\gamma > 1,$ the functions $\alpha(t)$ and
$\beta(t)$ are bounded, it follows from the expression for the total
energy. Really,
$$
{\mathcal E}(t)=(\alpha(t)^2+\beta(t)^2)G(t)+ {\mathcal E}_p(0)
G^{\gamma-1}(0) \frac{1}{G^{\gamma-1}(t)}\le {\mathcal E}(0), $$ it
implies
$$ \alpha^2(t)+\beta^2(t)\le  {{\mathcal E}(0)}{G_1(t)}-
{\mathcal E}_p(0) G_1^{1-\gamma}(0) G_1^{\gamma}(t)<+\infty, $$
$|{\rm div}\,{\bf V}|<\infty.$ Consequently, the density and
pressure are bounded for all solutions of the class we want to
construct.

If $\,\mu=0,\, $ the system (2.1.9) can be integrated and we can
obtain $$ \alpha(G_1)=\pm\sqrt{C_2 G_1^\gamma-C_1^2G_1^2+({\mathcal
E}(0)-lC_1)G_1-l^2/4}, $$ $$ \beta(t)=C_1G_1(t)+l/2, \qquad
-\int_{G_1(0)}^{G_1(t)} \frac{dG_1}{2G_1\alpha(G_1)} = t,
$$
where $C_1=\frac{2\beta(0)-l}{2G_1(0)},\quad C_2=(\alpha^2(0)+C_1^2
G_1^2(0)-({\mathcal E}(0)-lC_1)G_1(0)+l^2/4)/G_1^\gamma(0).$

In the cases $\, \mu>0 \,$ and $\,\mu=l=0\,, $ there is a unique
stable equilibrium in the origin. We have the following asymptotics of
the solution components as $t\to\infty:$

$$\mbox{if\,} \,\mu> 0, \,l=0,\,
\mbox{then}\,\alpha(t)\sim \frac{1}{2\gamma}t^{-1},\, \beta(t)\sim
C_1\left(\frac{\mu}{2K_1\gamma}\right)^{1/\gamma}t^{-1/\gamma}\exp\{-\mu
t\},$$$$ G_1(t)\sim \left(\frac{\mu}{2K_1\gamma}\right)^{1/\gamma}
t^{-1/\gamma};$$

$$\mbox{if\,} \,\mu>0,l\ne 0,\,\mbox{then}\,\alpha(t)\sim
\frac{1}{2\gamma}t^{-1},\, \beta(t)\sim
\frac{l}{2 \mu\gamma}t^{-1},\, G_1(t)\sim \left(
\frac{l^2+\mu^2}{2K_1\mu\gamma}\right)^{1/\gamma} t^{-1/\gamma};$$

$$\mbox{if\,}\, \mu=l=0,\,\mbox{then}\,
\alpha(t)\sim\left(\sqrt{G_1(0)/{\mathcal E}(0)}+t\right)^{-1},\,
G_1(t)\sim \alpha^2(t),\, \beta(t)= C_1 G_1(t),$$ where
$K_1=(\gamma-1){\mathcal E}_p(0) G_1^{1-\gamma}(0).$

Knowing $\alpha(t)$ and $\beta(t)$, in the standard
way  we can find components of density and entropy from the equations
(E1), (E3), linear with respect to $\rho$ and $p $, respectively.
We obtain that
$$ \rho(t,|r|,
\phi)=\exp(-2\int\limits_0^t \alpha(\tau) d\tau)
\rho_0(|r|\exp(-\int\limits_0^t \alpha(\tau) d\tau),
\phi+\int\limits_0^t \beta(\tau) d\tau), $$
$$
p(t,|r|, \phi)=\exp(-2\gamma\int\limits_0^t \alpha(\tau) d\tau)
p_0(|r|\exp(-\int\limits_0^t \alpha(\tau) d\tau),
\phi+\int\limits_0^t \beta(\tau) d\tau). $$
From (E2) and (2.1.9) we get that on the classical
solution of the  system (E1--E3) the relation
$${\bf
\nabla} p= - (\gamma-1)G_1^{1-\gamma}(0)E_p(0)G_1^\gamma(t)\rho {\bf
r}$$
must be satisfied.
Hence it follows that the components of the initial data
$\rho_0$ and $p_0$ must be asymmetrical and compatible,
i.e. connected as follows:
$${\bf
\nabla} p_0= - (\gamma-1)G_1(0)E_p(0)\rho_0 {\bf r}.\eqno(2.1.10)$$

For example, one can choose
$$p_0=\frac{1}{(1+|{\bf
r}|^2)^a},\,a=const>3,\eqno(2.1.11)$$
$$ \rho_0=
\frac{2a}{(\gamma-1)G_1(0)E_p(0)}\frac{1}{(1+|{\bf
r}|^2)^{a+1}}.\eqno(2.1.12)$$

\subsubsection{$A(t)$ of general form, ${\bf b}(t)=0$}

To consider the velocity field
(2.1) with the matrix
$$A(t)=\left(
\begin{array}{cc} a(t) & b(t)\\ c(t) &
d(t)\end{array}\right),
$$
we introduce the functionals
$$ G_x(t)=\frac{1}{2}\int\limits_{{\mathbb R}^2} \rho x^2
d{\bf x},\quad
G_y(t)=\frac{1}{2}\int\limits_{{\mathbb R}^2} \rho y^2 d{\bf x},
\quad G_{xy}(t)=\frac{1}{2}\int\limits_{{\mathbb R}^2} \rho xy
\,d{\bf x}.  $$
Then involve the auxiliary variables
$$G_1(t)={G_x(t)}{\Delta^{-(\gamma+1)/2}(t)},\,
G_2(t)={G_y(t)}{\Delta^{-(\gamma+1)/2}(t)},$$
$$
G_3(t)={G_{xy}}{\Delta^{-(\gamma+1)/2}(t)},\eqno(2.1.13)$$
where $\Delta(t)=G_x(t)G_y(t)-G_{xy}^2(t)$ is a positive
function on solutions to (E1--E3).  Note that the behaviour of
$\Delta(t)$ is governed by the equation
$$
\Delta'(t)=2(a(t)+d(t))\Delta(t),\eqno(2.1.14)$$
and the potential energy
$E_p(t)$ is connected with $\Delta(t)$ as follows:
$$E_p(t)=E_p(0)\Delta^{(\gamma-1)/2}(0)\Delta^{(-\gamma+1)/2}(t).$$

To find $G_1,\,G_2,\,G_3$, and  the elements of the matrix
$A(t)$ we get the system of equations
$$
G_1'(t)=((1-\gamma)a(t)-(1+\gamma)d(t))G_1(t)+2b(t)G_3(t), $$ $$
G_2'(t)=((1-\gamma)d(t)-(1+\gamma)a(t))G_2(t)+2c(t)G_3(t),
$$
$$
G_3'(t)=c(t)G_1(t)+b(t)G_2(t)-\gamma(a(t)+d(t))G_3(t),
$$
$$
a'(t)=-a^2(t)-b(t)c(t)+lc(t)-\mu a(t)+K_2G_2(t),\eqno(2.1.15)
$$
$$
b'(t)=-b(t)(a(t)+d(t))+ld(t)-\mu b(t)-K_2 G_3(t), $$
$$
c'(t)=-c(t)(a(t)+d(t))-la(t)-\mu c(t)-K_2G_3(t), $$ $$
d'(t)=-d^2(t)-b(t)c(t)-lb(t)-\mu d(t)+K_2G_1(t), $$
with
$K_2=\frac{\gamma-1}{2}
(E_p(0)\Delta^{(\gamma-1)/2}(0)).$

Now we obtain an analytical result on a behaviour of
the coefficients of the matrix $A(t)$ as $t\to\infty.$

\begin{utv 2.}
Suppose $\mu=l=0,$ that is ${\bf F=0}.$
If the system (E1--E3) has the solution  with the linear
profile of velocity ${\bf V}= A(t){\bf r},$ then
$A(t)\sim \frac{\delta}{t}I$
as $t\to\infty,$ with the identity matrix $I$ and a constant
$\delta>\frac{1}{2}.$
\end{utv 2.}

{\it Proof.} Go to the new variables
$a_1(t)=a(t)-d(t), b_1(t)= b(t)+c(t), c_1(t)=b(t)-c(t),
d_1(t)=a(t)+d(t), G_4(t)=G_1(t)+G_2(t), G_5(t)=G_1(t)-G_2(t).$
In  variables
$ a_1,b_1,c_1,d_1,G_3,G_4,G_5 $ system (2.1.15) has the form:
$$ a_1'(t)=-a_1(t)d_1(t)-K_2G_5(t),\eqno(2.1.16) $$ $$
b_1'(t)=-b_1(t)d_1(t)-2K_2G_3(t),\eqno(2.1.17)
$$
$$
c_1'(t)=-c_1(t)d_1(t),\eqno(2.1.18)
$$
$$
d_1'(t)=-\frac{1}{2}(a_1^2(t)+d_1^2(t))-\frac{1}{2}(b_1^2(t)-c_1^2(t))+
K_2G_4(t),\eqno(2.1.19)
$$
$$
G_3'(t)=-\gamma d_1(t)G_3(t)+\frac{1}{2}
b_1(t)G_4(t)-\frac{1}{2}c_1(t)G_5(t),\eqno(2.1.20) $$
$$
G_4'(t)=-\gamma
d_1(t)G_4(t)
+a_1(t)G_5(t)+2b_1(t)G_3(t),\eqno(2.1.21)$$
$$
G_5'(t)=-\gamma
d_1(t)G_5(t)+a_1(t))G_4(t)-2c_1(t)G_3(t).\eqno(2.1.22) $$

Let as $ t\to\infty $ the asymptotics of functions involved in the
system be the following:  $
a_1(t)\sim L_1 t^{l_1}, $ $b_1(t)\sim L_2 t^{l_2}, c_1(t)\sim L_3
t^{l_3}, $$\, d_1(t)\sim L_4 t^{l_4}, $$\,G_3(t)\sim N t^{q},\,
G_4(t)\sim M_1 t^{p_1}, $ $\,G_5(t)\sim M_2 t^{p_2},  $ where $L_i,
\,i=1,2,3;\, M_j,\, j=1,2; \, N $ are certain constants not
equal to zero.  Note that
$p_2\le p_1$ and in virtue of $\Delta>0$
the estimate $q\le p_1$ holds.

From (2.1.18) we get immediately that
$$l_3 L_3 t^{l_3-1}=-L_3 L_4 t^{l_3+l_4},$$
hence
$l_4=-1,\, L_4=-l_3.$

Taking the fact into account,
from (2.1.16) we have
$$l_1 L_1 t^{l_1-1}=-L_1 L_4 t^{l_1-1}-2K_2 M_2 t^{p_2}.$$
The following variants are possible:

$p_2<l_1-1, $ hence $l_1=-L_4,$ i.e. $l_1=l_3,$

$p_2=l_1-1, $ hence
$$l_1=-L_4-\frac{K_2M_2}{L_1}=l_3-\frac{K_2M_2}{L_1}.\eqno(2.1.23)$$

From (2.1.17) we have analogously
$$l_2 L_2 t^{l_2-1}=-L_2 L_4 t^{l_2-1}-2K_2 N t^q.$$
There are the variants:

$q<l_2-1, $ hence $l_2=-L_4,$ i.e. $l_2=l_3,$

$q=l_2-1, $ hence
$$l_2=-L_4-\frac{2K_2N}{L_2}=l_3-\frac{2K_2N}{L_2}.\eqno(2.1.24)$$

From (2.1.19) we get
$$-L_4 t^{-2}=-\frac{1}{2}(L_1^2 t^{2l_1}+L_4^2 t^{-2})-
\frac{1}{2}(L_2^2 t^{2l_2}-L_3^2 t^{2l_3})+
K_2M_1t^{p_1}.\eqno(2.1.25)$$

If $l_i<-1,\,i=1,2,3,$ then at $t\to\infty$
$$a(t)\sim\frac{1}{2}(L_1t^{l_1}+L_4 t^{-1}),\quad
b(t)\sim\frac{1}{2}(L_2t^{l_2}+L_3 t^{l_3}),$$
$$c(t)\sim\frac{1}{2}(L_2t^{l_2}-L_3 t^{l_3}),\quad
d(t)\sim\frac{1}{2}(L_4t^{-1}-L_1 t^{l_1}),$$
that is $A(t)\sim \frac{L_4}{2t}I,\, L_4=const>0.$

In the case $L_4=-l_3>1,$ therefore $\delta>\frac{1}{2}.$

We shall show below that others situations do not appear.

From (2.1.25) one can deduce the impossibility of the
situation $l_3=-1, l_1<-1,l_2<-1,p_1\le -2.$ Really, in the case
$l_3=-L_4=-1-\sqrt{1+L_3^2}<-1$ or
$l_3=-L_4=-1-\sqrt{1+L_3^2+2K_2M_1}<-1,$
it contradicts to the assumption.

Now let $l_3>-1.$ Then from (2.1.25)
$l_1=l_3$ or (and) $l_2=l_3.$
Suppose, for example, that $l_1=l_3.$
Then $p_2<l_1-1,$ from (2.1.22) it follows that
$$p_2M_2 t^{p_2-1}=-\gamma L_4 M_2 t^{p_2-1} + L_1M_1t^{p_1+l_1}
+2L_2Nt^{q+l_2},
\eqno(2.1.26)
$$
and therefore
$p_1+l_1=q+l_2,\,l_1\le\l_2,$
i.e. $l_2>-1.$
Besides,
$$L_1M_1=-2L_2N.\eqno(2.1.27)$$

From (2.1.20) we have
$$qN t^{q-1}=-\gamma L_4 N t^{q-1} + \frac{1}{2}L_2M_1t^{p_1+l_2}
-\frac{1}{2}L_3M_2t^{p_2+l_3}.
\eqno(2.1.28)
$$
As soon as $ q-1\le p_1-1, $ and  $ p_1-1 < p_1+l_2,$
then
$p_1+l_2=p_2+l_3\le p_1+l_3, $
and therefore
$l_2\le l_3,$
and as
$l_1=l_3, $ then $l_1=l_2=l_3$ and $p_1=q.$
Besides,
$$L_2 M_1=L_3M_2.\eqno(2.1.29)$$

Further, from (2.1.21) we obtain
$$p_1M_1 t^{p_1-1}=-\gamma L_4 M_1 t^{p_1-1} + L_1M_2t^{p_2+l_1}
-2L_3Nt^{q+l_3}.
\eqno(2.1.30)
$$
As $ q+l_3> p_1-1, $ then
$p_2+l_1=q+l_3,\, p_2=q, $
$$L_2 M_1=2L_3N.\eqno(2.1.31)$$

From (2.1.27) and (2.1.31) we have
$\frac{M_1}{M_2}=-\frac{L_2}{L_3},$ and from (2.1.29) we
obtain  $\frac{M_1}{M_2}=\frac{L_2}{L_3},$ i.e.
$L_2^2=-L_3^2,\, L_2=L_3=0$ in spite of the assumption.

Suppose now that $l_2=l_3.$ Then
$q<l_2-1,\,p_2\le l_1-1.$ From (2.1.20) we have
$\,q-1<l_2+p_1, $ since $l_2>-1,$ and therefore $l_2+p_1=l_3+p_2,\,
p_1=p_2,\,L_2M_1=L_3M_2.$

Therefore from (2.1.26)  we have
$p_1+l_1=q+l_2\le p_1+l_2,\,l_1\le l_2,\, L_2M_1=-2L_2N.$

From (2.1.30) we get
$p_2+l_1=q+l_3\le p_1+l_2,\,p_2=q,\,l_1=l_2,\, L_1M_2=-2L_3N.$
In that way,  we obtain the contradiction
analogous to the previous one.

Now let $l_3\le-1,$ а $l_1>-1 $ and (or)   $l_2>-1.$
Then from (2.1.25) it follows that
$p_1=2l_1 $ and (or)   $p_1=2l_2\, (p_1>-2). $
For example, if
$p_1=2l_1, $ then from (2.1.27) we get
$p_2-1<p_1+l_1,\, p_1+l_1=q+l_2\le p_1+l_2,\,l_1\le l_2.$
From (2.1.28)
$p_2+l_3\le p_1-1,\,l_2+p_1\ge l_1+p_1>p_1-1,$
therefore
$q-1=l_2+p_1.$
But
$q-1\le p_1-1,$
therefore $l_2\le-1,$ and $l_2\le-1$ in spite of the
assumption.

If $p_1=2l_2, $ then from (2.1.28) we get
$l_3+p_2\le p_1-1,$ therefore $p_1+l_2>=p_2+l_3,\,
q-1=p_1+l_2,\,l_1\le l_2.$ But $q-1\le p_1-1, $
therefore $l_2\le-1$ in spite of the supposition.

So, it remains the unique possibility:
$l_3\le -1,\,l_1\le -1 $ and (or) $l_2\le -1.$
In the case, as follows from (2.1.25)
$$l_3=-1\pm\sqrt{1-(\delta_1 L_1^2+\delta_2 L_2^2 -\delta_3 L_3^2
- 2\delta_4 K_2 M_1)},\eqno(2.1.32)$$
where $  \delta_i=1,$ if $
l_i=-1, $ and $ \delta_i=0 $ otherwise, $
i=1,2,3, \, $  $ \delta_4=1, $ if $p_2=-2 $  and $ \delta_4=0$
otherwise. Hence
$$\delta_1 L_1^2+\delta_2 L_2^2 -\delta_3 L_3^2
- 2\delta_4 K_2 M_1\le 1, \eqno(2.1.33)$$
if inequality (2.1.33) is strict, then $l_3<-1.$

We consider this case, that is
$l_3<-1,\,l_1=-1$ and (or)  $l_2=-1.$

If $l_1=-1\ne\l_3,$ then $p_2=l_1-1=-2.$
If $l_2=-1\ne\l_3,$ then
$q=l_2-1=-2.$  In that way, in any of these cases we have
$ p_1=-2.$

At $l_1=-1,\,l_2<-1,\,l_3<-1,$
from (2.1.26)
we obtain
$p_2=-\gamma L_4+\frac{L_1M_1}{M_2},$
$$\frac{L_1M_1}{M_2}=-(2+\gamma l_3),\eqno(2.1.34)$$
from (2.1.30)  we get
$p_1=-\gamma L_4+\frac{L_1M_2}{M_1},$
$$\frac{L_1M_2}{M_1}=-(2+\gamma l_3),\eqno(2.1.35),$$
and from (2.1.24)
$$
\frac{K_2M_2}{L_1}=l_3+1. \eqno(2.1.36)$$
From (2.1.34), (2.1.36) after excluding  $L_1$ and $M_2$
we obtain
$$K_2M_1=-(2+\gamma l_3)(l_3+1).\eqno(2.1.37)$$

As $K_2$ and $M_1$ are positive, and $l_3\le -1,$
then
$$l_3>-\frac{2}{\gamma}.\eqno(2.1.38)$$

At $\gamma>2$ it goes already to the contradiction.

Further, from (2.1.35) and (2.1.36) we have
$\frac{L_1^2}{K_2M_1}=-\frac{2+\gamma l_3}{l_3+1},$
this fact together with (2.1.37) gives
$$L_1^2= (2+\gamma l_3)^2. \eqno(2.1.39)$$

From (2.1.32), (2.1.37), (2.1.38) and (2.1.39) we obtain the
inequality $$-\frac{2}{\gamma}<-1-\sqrt{1-
(2+\gamma l_3)^2 -2 (2+\gamma l_3)(l_3+1)},\eqno(2.1.40)$$
which cannot be true at $\gamma>1.$

       If $l_1<-1,\,l_2=-1,\,l_3<-1,$ then from (2.1.26)
we get
$p_2=-\gamma L_4+\frac{2L_2N}{M_2},$
$$
\frac{2L_2N}{M_2}=-2-\gamma l_3, \eqno(2.1.41)
$$
from (2.1.28)
$q=-\gamma L_4+\frac{L_2M_1}{2N},$
$$
\frac{L_2M_1}{2N}=-2-\gamma l_3, \eqno(2.1.42)
$$
from (2.1.30)
$$p_2=-\gamma L_4=\gamma l_3.$$
In that way, $l_3=-\frac{2}{\gamma},$
this fact together with (2.1.41) or (2.1.42) contradicts to
the fact that $L_2, \,N$ and $M_1$ are not equal to zero.

    If $l_1<-1,\,l_2=-1,\,l_3=-1,$ then from
(2.1.26),(2.1.28) and (2.1.30) we get correspondingly
$$\frac{2L_2N}{M_2}+\frac{L_1M_1}{M_2}=-2-\gamma
l_3,\eqno(2.1.43)$$ $$ \frac{L_2M_1}{2N}=-2-\gamma l_3,
\eqno(2.1.44) $$ $$ \frac{L_1M_2}{M_1}=-2-\gamma l_3,
\eqno(2.1.45) $$ from (2.1.23), (2.1.24) $$
\frac{K_2M_2}{L_1}=\frac{2K_2N}{L_2}=l_3+1. \eqno(2.1.46)$$
In that way, multiplying (2.1.43) by (2.1.45), taking into
account (2.1.46) we get
$$ (2+\gamma
l_3)^2=L_1^2+\frac{M_2L_2^2}{M_1},$$ from (2.1.44), (2.1.45),
(2.1.46) $$M_1^2=M_2^2,\eqno(2.1.45)$$ that is
$$ (2+\gamma
l_3)^2=L_1^2\pm L_2^2.\eqno(2.1.48) $$ As above, from
(2.1.34), (2.1.37), (2.1.38) and (2.1.48) we get the
inequality (2.1.40) which cannot hold at $\gamma>1.$

It remains the unique possibility
$l_1=l_2=l_3=-1.$
In the case $p_2<-2,\,q<-2$ and, as follows from (2.1.26),
$p_1=p_2$ or $p_1=q,$ that is $p_1<-2.$ Therefore, as follows
from (2.1.32)
$$ L_1^2+L_2^2-L_3^2=1.\eqno(2.1.49)$$

If $p_1=p_2$, then $|M_2|<M_1.$
If  $p_1=p_2,\,q<p_1,$ then from (2.1.26), (2.1.30)
we obtain that
$p_2+\gamma=\frac{L_1M_1}{M_2},\,
p_1+\gamma=\frac{L_1M_2}{M_1},$
hence $M_1^2=M_2^2,$ it goes to the contradiction.

If $p_1=q,\,p_2<p_1,$
then
from (2.1.26), (2.1.28), (2.1.30)
we get
       $$L_1M_1=-2L_2N,\,
\lambda=
\frac{L_2M_1}{2N}=-\frac{2L_3N}{M_1},\eqno(2.1.50)$$
where
$\lambda= q+\gamma
=p_1+\gamma.$
It follows, in particular, that
$L_2^2=-\lambda L_1,\,L_1L_3=-\lambda L_2,
\,\lambda^2=-L_2L_3.$
Besides, if $G_1(t)\sim N_1
t^{q_1},\,G_2(t)\sim N_2 t^{q_2},\,t\to\infty, $ where $N_1, N_2
$ are some positive constants,
$q_1=q_2=p_1=q,\, N_1=N_2,\, M_1=2N_1, $ and $N_1^2\ge N^2 $ in
virtue of $G_1(t)G_2(t)-G_3^2(t)>0.$ In such way, from
(2.1.50) we get that $L_2^2\ge \lambda^2,\, \lambda^2 \ge
L_3^2,\,L_1^2\ge L_2^2,$ and taking into account (2.1.49),
$L_1^2\le 1,\, L_2^2\le 1,\,\lambda^2\le 1.$ That is if
$\lambda>1\,(\gamma>3),$ the conditions mentioned in  the
paragraph cannot hold together.

At last, consider the case $p_1=p_2=q.$ Then from
(2.1.26), (2.1.28), (2.1.30) we have
$\lambda=\frac{L_1M_1}{M_2}+\frac{2L_2N}{M_2}=\frac{L_2M_1}{2N}-
\frac{L_3M_2}{2N}=\frac{L_1M_2}{M_1}-\frac{2L_3N}{M_1}, $
that is the system of linear homogeneous equations with respect to
the variables $M_1,\,M_2,\,N$ $$ L_1M_1-\lambda M_2+2L_2N=0,$$ $$
L_2M_1-L_3M_2-2\lambda N=0,$$ $$ \lambda M_1-L_1M_2+2L_3N=0.$$
For the existence of its nontrivial solution the determinant of
the system must be equal to zero, i.e.
$$
\lambda L_1^2-2\lambda L_2L_3+L_1L_2^2+L_1L_3^2-\lambda^3=0.$$
Taking into account (2.1.49), involve the function with
respect to the variables $L_2$ and $L_3$, where $\lambda$ plays the
role of parameter, namely
$$\Psi_{\lambda}(L_2,L_3)=\lambda(1-L_2^2+L_3^2)-2\lambda L_2 L_3 +
(L_2^2+L_3^2)\sqrt{1-L_2^2+L_3^2}-\lambda^3.$$
By the standard methods one can show
that the function is not equal to zero at
$\lambda>1,$ i.e. at $\gamma\ge 3$
(remember, that $p_1<-2$).

So, it remains to investigate the case
$1<\gamma\le 3.$

Consider equation (2.1.14), which can be written as
$$ \Delta'(t)=2d_1(t)\Delta(t).$$
If we suppose that
$ \Delta(t)\sim const\cdot t^m,\,t\to\infty,\,m $ is a
constant, then $ m=2L_4=2.$ From (2.1.14) we can get
$E_p(t)\sim const\cdot t^{1-\gamma}(\to 0), \,t\to\infty.$
Therefore,
if we denote $E$ the quantity of the total energy of the
system, then $E_k(t)=E-E_p(t)\sim
E(1-C_1 t^{1-\gamma}),$ here and further $C_i $
are
some positive constants.  But
$E_k(t)\sim C_2 G_k(t) t^{-2}\Delta^{\frac{\gamma+1}{2}}\sim C_3 G_k
t^{\gamma-1},$ where $G_k  $ is at least one of functions
$G_1,\,G_2$ or $G_3.$ That is $G_k\sim
C_4(1-C_1 t^{1-\gamma})t^{1-\gamma}\sim C_5 t^{1-\gamma}.$ But then
$1-\gamma\le p_1<-2,\, \gamma>3, $ and $\gamma$
does not belong to the interval under consideration.

So the proof of the proposition is over.

\medskip

{\sc Remark 2.1.1} One can  show more shortly, that in the physical
case $1<\gamma \le 2$ the situation $l_1=l_2=l_3=-1$ if impossible.
Taking into account (2.1.13) we have
$G_1G_2-G_3^2=\Delta^{-\gamma}\sim const\cdot t^{-2\gamma},\,
t\to\infty.$ But the degree of the leading term of the expression
$G_1G_2-G_3^2$ is not greater then $2p_1,$ therefore $p_1\ge-\gamma,$
as $p_1<-2,$ then $\gamma>2.$
\medskip

{\sc Remark 2.1.2} Actually $l_3=-L_4=-2,\,
G_1G_2-G_3^2\sim const\cdot t^{-4\gamma},\,t\to\infty.$

\medskip

\subsubsection{$A(t)$ of special form, ${\bf b}(t) \ne 0$}

Suppose that the velocity  field is of form (2.1.8).

In the case we need to analyze rather complicated system
of 7 equations:

$$G'(t)=2\alpha(t)G(t)+b_1(t)N_1(t)+b_2(t)N_2(t),$$

$$N_1'(t)=\alpha(t)N_1(t)+\beta(t)N_2(t)+b_1(t){\mathcal M},$$

$$N_2'(t)=\alpha(t)N_2(t)-\beta(t)N_1(t)+b_2(t){\mathcal M},$$

$$(\alpha(t)N_1(t)+\beta(t)N_2(t)+b_1(t){\mathcal M})'=$$
$$\quad-\mu
(\alpha(t)N_1(t)+\beta(t)N_2(t)+b_1(t){\mathcal M})-l
(\alpha(t)N_2(t)-\beta(t)N_1(t)+b_2(t){\mathcal M}),$$

$$(\alpha(t)N_2(t)-\beta(t)N_1(t)+b_2(t){\mathcal M})'=
$$$$\quad-\mu
(\alpha(t)N_2(t)-\beta(t)N_1(t)+b_2(t){\mathcal M})+l
(\alpha(t)N_1(t)+\beta(t)N_2(t)+b_1(t){\mathcal M}),$$

$$(b_1(t)N_1(t)+b_2(t)N_2(t)+2\alpha(t)G(t))'=$$
$$\quad-\mu
(2\alpha(t)G(t)+b_1(t)N_1(t)+b_2(t)N_2(t))+l
(2\beta(t)G(t)+b_2(t)N_1(t)-b_1(t)N_2(t))+
$$$$\quad (b_1^2(t)+b_2^2(t)){\mathcal M}+
2(\gamma-1)KG^{1-\gamma}(t)+2G(t)(\alpha^2(t)+\beta^2(t))+$$$$\quad
2\alpha(t)b_1(t)N_1(t)+2\beta(t)b_1(t)N_2(t)+2\alpha(t)b_2(t)N_2(t)
-2\beta(t) b_2(t)N_2(t),$$

$$(b_2(t)N_1(t)-b_1(t)N_2(t)+2\beta(t)G(t))'=$$$$\quad
-l(2\alpha(t)G(t)+b_1(t)N_1(t)+b_2(t)N_2(t))
-\mu(2\beta(t)G(t)+b_2(t)N_1(t)-b_1(t)N_2(t)).$$

We could write the normal form of the  system, but it is very long.
Mention only  that we can always express the derivatives of $\alpha,
\beta, b_1, b_2$ explicitly, because the determinant of the
corresponding algebraic system is equal to $$4G^2(t){\mathcal
M}^2-(N_1(t)^2 +N_2(t)^2)^2.$$ As the H$\rm \ddot o$lder inequality
shows, it is positive for the solutions we consider.

The compatibility condition in the case has the form
$${\bf \nabla} p_0=\rho_0(c_1{\bf r}+c_2{\bf r}_\bot +{\bf c}_0)$$
with scalar constants $c_1,\,c_2$ and constant vector ${\bf c}_0.$
\medskip

\subsubsection{$A(t)$ of general form, ${\bf b}(t) \ne 0$}

The corresponding system consists of 11 equations for
unknown functions $a(t),$ $ b(t),$ $ c(t),$ $d(t),$ $b_1(t),$
$b_2(t),$ $ G_x(t),$ $ G_y(t),$ $G_{xy}(t),$ $ N_1(t),$ $ N_2(t).$ It
seems that it can be analyzed only numerically.

As soon as the components of velocity field are found,
we can always find the density and the pressure.
Namely,
$$\rho(t,{\bf x})=\exp (-\int\limits_0^t{\rm tr} A(\tau) d\tau)
\rho_0(M^{-1}(t){\bf x}-\int\limits_0^t M^{-1}(\tau){\bf
b}(\tau)d\tau),\eqno(2.1.51)$$
$$p(t,{\bf x})=\exp (-\gamma\int\limits_0^t{\rm tr}
A(\tau) d\tau) p_0(M^{-1}(t){\bf x}-\int\limits_0^t M^{-1}(\tau){\bf
b}(\tau)d\tau),\eqno(2.1.52)$$
where $M(t)$ is a  solution
to the matrix equation $\displaystyle \frac{d X }{dt}=A(t){X}(t)$
\cite{Bellman}.

The compatibility condition in the case has the form
$${\bf \nabla} p_0({\bf x})=\rho_0({\bf x})(C{\bf x}+{\bf c_0})$$
with a constant matrix $C$ and a constant vector
${\bf c_0}.$

\medskip

\subsection{The system in the physical space ($n=3$)}

Now we consider the right-hand side of (E1), corresponding to
the Coriolis force and the Rayleigh friction in the real physical
space.  Namely,  $${\bf f}=-\mu{\bf
V}+\delta{[\bf V\times{\bf\omega}]}=
(-\mu I+\delta\left(\begin{array}{ccc}
0&\omega_3&-\omega_2\\-\omega_3&0&\omega_1\\
\omega_2&-\omega_1&0
\end{array}\right)){\bf V},
\eqno(2.2.1)$$
where $\delta=0$ or $1$, $\bf \omega$ is a
constant vector $(\omega_1,\omega_2, \omega_3),$ $\mu\ge 0$ is a
constant, $[\bf a\times b]$ is the cross product of the vectors $\bf
a$ and $\bf b$.

Besides of $G(t)=\frac{1}{2}\int
\limits_{{\mathbb R}^3}\rho |{\bf r}|^2\,d{\bf x} $ and
$F_1(t)=\int\limits_{{\mathbb R}^3}({\bf V, r})\rho\,d{\bf x},$
the
generalization of the functionals, involved in 2.1.1, to the 3D case,
we consider
$$ \tilde F_2(t)=\int\limits_{{\mathbb R}^3}({\bf V}, {[\bf
\omega \times r]})\rho\,d{\bf x},$$ $$ F_3(t)=\int\limits_{{\mathbb
R}^3}({[\bf V\times \omega]}, {[\bf \omega \times r])}\rho\,d{\bf
x},$$ $$H(t)=\frac{1}{2}\int\limits_{{\mathbb R}^3}({[\bf \omega
\times r]})^2\rho\,d{\bf x}.$$
For the smooth solutions to the system (E1--E3) with the right-hand
side (2.2.1) we have
$$G'(t)=F_1(t),$$$$F_1'(t)=2 {\mathcal E}_k(t)+3(\gamma-1){\mathcal E}_p(t)+
\delta \tilde F_2(t)-\mu F_1(t),$$
$$\tilde F'_2(t)=\delta F_3(t)-\mu\tilde F_2(t), $$
$$H'(t)=0.$$
We seek the solution with a special linear velocity profile
$${\bf V}= \alpha(t){\bf r}+ \beta(t){[\bf r \times \omega]}
=(\alpha(t) E+\beta(t)\left(\begin{array}{ccc}
0&\omega_3&-\omega_2\\-\omega_3&0&\omega_1\\
\omega_2&-\omega_1&0
\end{array}\right)){\bf r}.$$
In this case $$F_1(t)=2\alpha(t)G(t),$$$$\tilde F_2(t)=-2\beta(t) H(0),
$$$$F_3(t)=-2\alpha H(0),$$
$${\mathcal E}_k(t)=\alpha^2(t)G(t)+\beta^2(t)H(0),$$
$$ {\mathcal
E}_p(t)=const\cdot G^{-\frac{3(\gamma-1)}{2}}(t).$$
Changing $G^{-1}(t)=G_1(t)$, we obtain the closed system
of ODE
$$\alpha'(t)=-\alpha^2(t)-\mu\alpha(t)+\beta^2(t)G_1(t)H(0)-\delta
\beta (t) G_1(t) H(0)+\frac{3(\gamma-1)}{2}\tilde K
G_1^{\frac{3\gamma-1}{2}}(t),$$
$$\beta'(t)=\delta\alpha(t)-\mu\beta(t),\quad
G_1'(t)=-2\alpha(t)G_1(t), $$ with $\tilde K= {\mathcal
E}_p(0)G^{\frac{3(\gamma-1)}{2}}(0).$ Since the total energy is
nonincreasing,  we can draw similarly to the 2D case that
$|\alpha(t)|<\infty, |\beta(t)|<\infty.$

In
the situation with $\mu=\delta=0,\,$ with $\mu>0,\delta=0$ and
with $\mu>0,\delta=1,$ there is a stable equilibrium in the origin
$\alpha(t)=\beta(t)=G_1(t)=0.$

If $\mu=0,\delta=0,$ we have

$$
\alpha(t)\sim t^{-1}, \beta(t)=\beta(0)=const, G_1(t)\sim const\cdot
t^{-2}, t\to\infty.$$

If $\mu>0,\delta=0,$ then
$$ \alpha(t)\sim \frac{1}{3\gamma-1} t^{-1},
\beta(t)=\beta(0)\exp(-\mu t), G_1(t)\sim {const}\cdot
t^{-\frac{2}{3\gamma-1}}, t\to\infty.$$

In the case  $\mu>0,\delta=1,$
$$ \alpha(t)\sim \frac{1}{3\gamma-1} t^{-1},
\beta(t)\sim \frac{1}{\mu(3\gamma-1)} t^{-1},
G_1(t)\sim {const}\cdot
{t^{-\frac{2}{3\gamma-1}}}, t\to\infty.$$

What about other component of solution,
they can be found by the
same formulas, as in the 2D case ((2.1.50 -- 2.1.51)), and must
satisfy the compatibility condition
$${\bf
\nabla} p_0= - (\gamma-1)G_1(0)E_p(0)\rho_0 {\bf r},$$
the same as (2.1.10).

{\sc Remark 2.2.1} We can also consider the velocity with linear
profile of general form, but the system of integral functionals is
complicated and it is difficult to analyze it.

\subsection{Theorem on the interior solutions}

Further we need to obtain the symmetric form  of system (E1-E3).
For this purpose we consider the entropy $S,$ connected with the
components of solution by the state equation
$p = e^S \rho^\gamma.$

Thus, instead of (E1-E3) we obtain other system
$$ \partial_t \rho + {\rm div}\, (\rho{\bf V})=0, \eqno (2.3.1) $$ $$
\rho\partial_t {\bf V}+(\rho{\bf V},{\bf \nabla})\,{\bf
V}+{\bf\nabla} p = \rho {\bf f}({\bf x}, t, {\bf V}, \rho, S), \eqno
(2.3.2)$$
$$ \qquad \partial_t S +({\bf V},{\bf\nabla} S)=0. \eqno (2.3.3)$$

Here $\rho(t,{\bf x}) ,\,{\bf V}(t,{\bf x}),\, S(t,{\bf x})$ are
components of the solution, given in ${\mathbb R}_+\times {\mathbb
R}^n,\, n\ge 1$ (density, velocity and entropy,
respectively).

The systems (E1--E3) and (2.3.1--2.3.3) are equivalent for $\rho>0.$

Put the Cauchy problem for (2.3.1--2.3.3):
$$ \rho(0,{\bf x})=\rho_0({\bf x})\ge 0,
\, {\bf V}(0,{\bf x})={\bf V}_0({\bf x}),\, S(0, {\bf x})=S_0({\bf
x}).\eqno (2.3.4)$$

Denote  $p_0({\bf x}) = e^{S_0({\bf x})} \rho_0^\gamma({\bf x}),$
and remark that the first of the Cauchy conditions can be replaced
to $ p(0,{\bf x})=p_0({\bf x})\ge 0.$

\vskip1cm
%\begin
{\bf Definition.}
{\it
We shall call the global-in-time classical solution
$(\bar \rho(t,{\bf x}),\,$
$ \bar{\bf V}(t,{\bf x}), \,$
$\bar S(t,{\bf
x}))$ to the system (2.3.1--2.3.3)  {\bf  the interior
solution},
if $$(\bar p^{\frac{\gamma-1}{2\gamma}}, \,\bar {\bf
V},\,\bar S)\in C^1({\mathbb R}_+\times {\mathbb R}^n),$$
where
$\bar p = e^{\bar S} \bar\rho^\gamma$, and any solution $(\rho(t,{\bf
x}),\, {\bf V}(t,{\bf x}),\,S(t,{\bf x}))$ to the Cauchy problem
(2.3.1--2.3.4) with the sufficiently small norm
$$\|(p_0^{(\gamma-1)/2\gamma}({\bf x})- \bar
p^{(\gamma-1)/2\gamma}(0,{\bf x}),\,{\bf V}_0({\bf x})-\bar
{\bf V}(0, {\bf x}),\, S_0({\bf x})-\bar S(0,{\bf x})
\|_{H^m({\mathbb R}^n)},$$ $m>n/2+1,$ is also smooth,
such that $(p^{\frac{\gamma-1}{2\gamma}}, \,{\bf V},\, S )\in
C^1({\mathbb R}_+\times {\mathbb R}^n).$}

\vskip1cm
Note that the trivial solution is not interior at least for $\bf F=f
=0$.
But a set of the interior solutions is not empty. In
\cite{Grassin} (a generalization of
\cite{Serre}) for ${\bf f}={\bf 0}$ it was shown that the solution
$(0,\, \bar {\bf V}(t,{\bf x}),\,const)$ to (2.3.1--2.3.4) is
interior, provided $\bar {\bf V}(t,{\bf x}) $ is a globally smooth
solution to the equation $\partial_t {\bf V}+({\bf V},{\bf
\nabla})\,{\bf V} = 0$ such that the spectrum of its Jacobian is
separated  initially from the real negative semi-axis, $D\bar{\bf
V}(0,{\bf x})\in L^\infty({\mathbb R}^n),\, D^2\bar{\bf V}(0,{\bf
x})\in H^{m-1}({\mathbb R}^n)$ (we denote by $D^k$ the vector of all
spatial derivatives of the order $k$).

The result is clear from a physical point of view: the velocity
field with a positive divergency spreads the initially
concentrated small mass, that  prevents the singularity formation.

As we show below, some of solutions with linear profile of
velocity are interior
with the density not close to zero.

\vskip1cm

For the solutions to
(2.3.1--2.3.3) with a finite moment we have $\inf \rho=\inf p= 0,$
thus we need to use the symmetrization proposed in \cite{Makino}.
After involving the new variable $\displaystyle \Pi=\kappa
p^{\frac{\gamma-1}{2\gamma}},\,\kappa=\frac{2\sqrt{\gamma}}{\gamma-1},$
we obtain the symmetric form of the system (2.3.1--2.3.3) with the
solution ${\bf U}= (\Pi,\,{\bf V},\,S)$:
$$ \exp(\frac{S}{\gamma })(\partial_t +({\bf V},{\bf
\nabla}))\Pi + \frac {\gamma -1}{2}\exp(\frac{S}{\gamma})\Pi div\,
{\bf V}=0,\eqno(2.3.5) $$ $$ (\partial_t +({\bf V},{\bf
\nabla})){\bf V}+\frac {\gamma -1}{2}\exp
(\frac{S}{\gamma})\Pi{\bf\nabla} \Pi = {\bf f}_1( t, {\bf x},\Pi,
{\bf V}, S),\eqno(2.3.6) $$ $$ (\partial_t +({\bf V},{\bf \nabla}))
S=0,\eqno(2.3.7) $$ where ${\bf f}_1={\bf f}( t, {\bf
x},e^{-\frac{S}{\gamma}}\left(\frac{\Pi}{\kappa}\right)^{\frac{2}{\gamma-1}},
{\bf V}, S).$
Let  $\Pi_0:=\kappa
p_0^{\frac{\gamma-1}{2\gamma}},$
 ${\bf U}_0:=(\Pi_0,{\bf
V_0}, S_0).$
In this way, the Cauchy problem (2.3.1--2.3.4) is transformed to
the  problem (2.3.5--2.3.7), with the initial data
${\bf U}(0,{\bf
x})={\bf U}_0({\bf x}).$
Further,
suppose that the system (2.3.5--2.3.7) has a classical solution
$\bar{\bf U}:=
(\bar\Pi,\,\bar{\bf V},\, \bar S),$
where
$\bar{\bf V}:= A(t){\bf r}+{\bf b}(t),$ $\bar\Pi:=\kappa
\bar p^{\frac{\gamma-1}{2\gamma}}.$
Denote also
${\bf U}-\bar{\bf U}:={\bf u}:=
(\pi,\,{\bf v},\,s).$

Before formulating Theorem 2.1 we perform certain
transformations. Firstly we go on to the system with the solution
${\bf u}.$ Then,
following \cite{Serre}, we carry out
the nondegenerate change of variables such that the infinite
semi-axis of time turns to semi-interval.

Assume that there exists a nondegenerate square $(n\times n)$ matrix
$A_1(t)$ such that $A_1(t)(A_1^{-1}(t))'=A(t)$, moreover,
$A_1(t)A(t)=A(t)A_1(t)$. Choose a positive decreasing function $
\lambda(t)$ such that the integral $\int_0^\infty \lambda(\tau)
d\tau$ converges to the finite value $\sigma_\infty,$ and set
$\sigma(t)=\int_0^{t} \lambda(\tau) d\tau.$ Let $(\sigma,{\bf y}):=
(\sigma (t), A_1(t){\bf x})$ be new variables.  Then
${\bf\nabla}_{\bf x}=A_1^*{\bf\nabla}_{\bf y},\,{\rm div}_{\bf
x}{\bf V}= {\rm div}_{\bf y} A_1{\bf V},
\,\partial_{\sigma}=\lambda^{-1}(t)(\partial_t+A{\bf
r}{\bf\nabla}_{\bf x})$.  In that way, the semi-infinite axis of
time goes to the semi-interval $[0,\sigma_\infty).$ Note that there
exists the inverse function $t=t(\sigma).$  Further, involve the
variables $ {\bf W}=\lambda^{-1}(t) A_1(t){\bf v},\,
P=\lambda^{-q}(t)\pi, \, $ the constant $q$ will be defined below.
Let ${\mathcal U}(\sigma,{\bf y}):=(P,{\bf W},s),\,\bar P:=
\lambda^{-q}\bar\Pi,\, \bar {\bf W}:=\lambda^{-1}A_1\bar{\bf V},\,
\bar{\mathcal U}:=(\bar P, \bar W, \bar S).$

So we get a  system
$$ (\partial_\sigma +({\bf W},{\bf \nabla_y}))P +
\frac {\gamma-1}{2}(P+\bar P) {\rm div}_{\bf y}\, {\bf W}
=
$$
$$
-({\bf W},{\bf \nabla_y})\bar P
-
PQ_1-\lambda^{-1}(A_1{\bf b},{\bf \nabla_y}
P),\eqno(2.3.8) $$ $$ (\partial_\sigma +({\bf W},{\bf
\nabla_y})){\bf W}+ \frac {\gamma -1}{2}\Psi(S,\sigma) (P+\bar
P){\bf\nabla_y} P = $$ $$ - \frac {\gamma -1}{2} \Psi(S,\sigma)
(P+\bar P){\bf\nabla_y} \bar P
+
\frac {\gamma -1}{2}
\Psi(\bar S,\sigma)
\bar P{\bf\nabla_y} \bar P-
$$$$
-Q_2{\bf
W}+G-\lambda^{-1}
(A_1{\bf b},{\bf \nabla_y}){\bf W}
,\eqno(2.3.9) $$ $$ (\partial_\sigma +({\bf W},{\bf
\nabla_y})) s =-({\bf W},{\bf \nabla_y})\bar S-
\lambda^{-1}(A_1{\bf b},{\bf \nabla_y}
s).\eqno(2.3.10) $$ Here we
take into account that $t=t(\sigma),\, {\bf
x}=A_1^{-1}(t(\sigma)){\bf y},\,$ $\lambda=\lambda(t(\sigma)),\,$$
A=A(t(\sigma)),\, $$ A_1=A_1(t(\sigma)),\, $${\bf b}={\bf
b}(t(\sigma))$ and denote
$$
Q_1=
Q_1(t(\sigma)):=\lambda(t)^{-1}(\frac{\gamma-1}{2} {\rm tr} A(t)+
q (\ln \lambda(t))'),$$
$$Q_2=Q_2(t(\sigma)):=\lambda(t)^{-1}((\ln \lambda(t))' E
+A(t)+A_1(t)A(t)A_1^{-1}(t)),$$
$$\Psi(S,\sigma):=
\exp
(\frac{S}{\gamma})R A_2,$$
$$A_2=A_2(t(\sigma)):=
A_1(t)A_1^*(t)({\rm det }A_1(t))^{-2/n}$$
$$R=R(t(\sigma)):= \lambda^{2q-2}(t) ({\rm det}\,A_1(t))^{2/n},$$
$$G=G(\sigma,{\bf y},{\mathcal U},\bar {\mathcal U})
:=\lambda^{-2}(t)A_1(t)({\bf f}_1( t, {\bf x},
\lambda^{q}(t)(\bar P+ P),$$$$\lambda(t)A_1^{-1}(t) (\bar{\bf W}+{\bf
W}), (\bar S+s))-{\bf f}_1( t,{\bf x},
\lambda^{q}(t)\bar P, \lambda(t)A_1^{-1}(t)\bar{\bf W}, \bar S)).$$
Note that  $A_2$ is a family  of invertible matrices.

The initial data ${\bf U}_0$ and the vector-function $\bar {\bf U}$
induce the Cauchy data ${\mathcal U}_0({\bf y})$ for the system
(2.3.8--2.3.10) as follows: ${\mathcal U}_0({\bf
y})=\lambda^{-q}(\Pi_0({A_1^{-1}\bf y})- $ $\bar\Pi(0,A_1^{-1}{\bf
y} ),$ $\lambda^{-1}A_1({\bf V}_0(A_1^{-1}{\bf y})-$ $\bar{\bf
V}(0,A_1^{-1}{\bf y})),$ $ S_0(A_1^{-1}{\bf y})- \bar
S(0,A_1^{-1}{\bf y}))|_{\sigma=0}.$

After multiplying (2.3.9) by $\Psi^{-1}(S,\sigma)$ the system
becomes symmetric hyperbolic. We can apply the theorem on a local
existence of the unique solution to the Cauchy problem for symmetric
hyperbolic systems \cite{Kato} to the problem (2.3.8--2.3.10),
${\mathcal U}(0,{\bf y})={\mathcal U}_0({\bf y})$, provided
$\Psi^{-1}(S,\sigma)$ is uniformly positive. Sometimes
$\Psi^{-1}(S,\sigma)$ is really uniformly positive (see
\cite{Serre}), and in that cases we can conclude that if ${\mathcal
U}_0\in H^m({\mathbb R}^n), m>1+n/2,$ then ${\mathcal U}\in \cap
_{j=0}^m C^j([0,\sigma_*);H^{m-j}({\mathbb R}^n)), \sigma_*>0.$
Since $\sigma_*\to \infty,$ as $\|{\mathcal U}_0 \|_{H^m}\to 0,$
then choosing the initial data ${\mathcal U}_0$ small in the Sobolev
$H^m -$ norm, we can extend the time of existence of the smooth
solution $\mathcal U$ to the Cauchy  problem up to $\sigma_*\ge
\sigma_\infty.$ In this case we conclude that ${\mathcal U}\in \cap
_{j=0}^m C^j([0,\sigma_\infty);H^{m-j}({\mathbb R}^n))$ and ${\bf
u}\in \cap _{j=0}^m C^j([0,\infty);H^{m-j}({\mathbb R}^n)).$

However, generally speaking, $\Psi^{-1}(S,\sigma)\to 0$ as $t\to
\infty,$ therefore we can extend the time of existence of the smooth
solution $\mathcal U$ to the Cauchy  problem up to $[0,\sigma_*),$
for any $\sigma_*<\sigma_\infty$, choosing the initial data small in
the Sobolev $H^m -$ norm, but we cannot guarantee the existence of
the solution for $[0,\sigma_\infty).$ Thus, we need  to act
differently.

If $t_*$ is the supremum of the time of existence of the smooth
solution ${\bf U}$ for the Cauchy problem (2.3.5--2.3.7), ${\bf
U}(0,{\bf x})={\bf U}_0({\bf x})$ (the solution ${\bf u}$ preserves
the smoothness during the same time), then $\sigma_*=\sigma(t_*).$
If $\sigma_*<\sigma_\infty (t_*<\infty),$ then $\displaystyle\lim
\sup\limits_{\sigma\to\sigma_*}(\|{\mathcal U}\|_{L_\infty}(\sigma)+
\|{\bf\nabla_y}{\mathcal U}\|_{L_\infty}(\sigma))=+\infty $
$(\displaystyle\lim \sup\limits_{t \to t_*}(\|{\bf
u}\|_{L_\infty}(t)+ \|{\bf\nabla_x}{\bf
u}\|_{L_\infty}(t))=\infty)$. Therefore, if we prove that
$\|{\mathcal U}\|_{H^m}(\sigma)<+\infty, \sigma\in [0,
\sigma_\infty)$, then we establish that $\sigma_*=\sigma_\infty
(t_*=\infty).$

  Let
$$D^J=\left(\frac{\partial}{\partial y_1}\right)^{j_1}...
\left(\frac{\partial}{\partial y_n}\right)^{j_n},\,
\sum\limits_i j_i=|J|.$$
Further,  ${\bf f}_{1\bf V}$ stands for the
matrix $\|\frac{\partial f_{1i}}{\partial V_j}\|$, $\nabla_\Pi{\bf
f}_1$ and $\nabla_S{\bf f}_1$ for the derivatives of $\bf f_1$ with
respect to $\Pi$ and $S$, correspondingly.

\newtheorem{theorem 2.}{Theorem 2.}
\begin{theorem  2.}
Let the function ${\bf
f}_1( t, {\bf x},\Pi, {\bf V}, S)$
have the derivatives with respect to all arguments up to the
order $m+1,
m>n/2+1$.    Suppose
that the system (2.3.1--2.3.3) has a global in time classical solution
$(\bar\rho, \bar{\bf V}, \bar S)$ with linear profile of velocity
$\bar{\bf V}=A(t){\bf r}+{\bf b}(t)$ such that

a) $\bar p^{\frac{\gamma-1}{2\gamma}}(t,{\bf x})\in
\cap _{j=0}^{m+1} C^j([0,\infty);H^{m-j}({\mathbb R}^n));$

\quad $D\,\bar S(t,{\bf x})
\in \cap _{j=0}^{m} C^j([0,\infty);H^{m-j}({\mathbb R}^n));$

b)there exists a matrix $A_1(t)$ such that
$A(t)=A_1(t)(A_1^{-1}(t))', \,$$A(t)A_1(t)=A_1(t)A(t),$ $\, \det
A_1(t)>0,$ for $t\ge t_0 \ge 0.$

Let us assume that there exist a smooth real-valued decreasing nonnegative
function
 $\lambda(t)$,
a constant $q$ and a
matrix $U_\phi(t)$ with real coefficients, having the following
properties:
$$\int_{t_0}^{+\infty} \lambda(\tau)
d\tau<\infty,\eqno(2.3.11)$$
$$
\quad\int_{t_0}^{+\infty} \lambda^q(\tau) (\det
A_1(\tau))^{1/n}\, d\tau<\infty,\quad \eqno (2.3.12)$$
$$R(t)'\ge 0 \quad \mbox{and}\quad Q_1(t)\ge 0 \quad\mbox{for} \quad
t\ge
t_0\ge 0,\eqno(2.3.13)$$
$$
A_2^{-1}(t) U_\phi(t)\mbox{ is a skew-symmetric matrix,}
$$
the following functions, vector-functions and matrices
are bounded in $[t_0,\infty)\times {\mathbb R}^n$ (under bounded
$(\Pi, {\bf V}, S))$:  $$\lambda^{-q}(t)\exp(-\frac{\gamma-1}{2}{\rm
tr}A(t)), \qquad\exp(-\int_{t_0}^t
A(\tau)d\tau),\eqno(2.3.14)$$
$$\lambda^{-q}(t) {\bf b}(t),
\qquad\lambda^{-1}(t) A_2^{-1}(t)A_2't)\eqno(2.3.15) $$$$
\lambda^{-1}D^J(\nabla_\Pi{\bf f}_1),\quad
\lambda^{-(q+1)}D^J(\nabla_S{\bf
f}_1),\eqno(2.3.16)$$
$$Q_3:=(Q_2(t)-\lambda^{-1}(t)D^J(A_1(t)\nabla_{\bf V}{\bf f}_1)
A_1^{-1}(t)- U_\phi(t)), \,|J|=0,1,\dots,m.  $$

Then the solution $(\bar\rho, \bar{\bf V}, \bar S)$ is interior.
Moreover,
$$(p^{(\gamma-1)/2\gamma}-
\bar p^{(\gamma-1)/2\gamma}, \, {\bf V}-\bar {\bf V},\,S-\bar S
)\in \cap _{j=0}^m C^j([0,\infty);H^{m-j}({\mathbb R}^n)).$$
\end{theorem 2.}

\vskip1cm

{\sc Remark 2.3.1.}
{\it If $Q_3$ depends only on time, then we can require only the
nonnegativity of this expression for $t\ge t_0.$ }

{\it
Proof.}
Let $\sigma\in [0,\sigma_*)$. Suppose for the sake of simplicity that
$t_0=0,$ otherwise we can consider the Cauchy problem at $
\sigma_0=\sigma(t_0)<\sigma_*$ with the initial data induced by the
 Cauchy data ${\mathcal U}_0$. Note that choosing ${\mathcal U}_0$
sufficiently small in the $H^m $ - norm we can prolong the smooth
solution ${\mathcal U}$ up to $\sigma_0$. Denote  $(X,\,{\bf
Y},\,Z):  R^n\to R\times R^n \times R $ a vector-function from
$L_2({\mathbb R}^{n+2})$. Introduce a norm $$ [X,{\bf
Y},Z]^2(\sigma):=\int\limits_{{\mathbb R}^n}((X^2+Z^2)+ {\bf
Y}^*\Psi^{-1}(S,\sigma){\bf Y}) d{\bf y}.$$

It is equivalent to the usual $L^2({\mathbb R}^{n+2})$ - norm, but
not uniformly in time. For any $p\in N$ we define
$$E_p(\sigma)=\frac{1}{2}\sum\limits_{|J|=p} [D^J
P,\,D^J {\bf W},\, D^J s]^2(\sigma),\qquad
F_m(\sigma)=\sum\limits_{p=0}^m E_p(\sigma).$$

Let us compute
$$\frac{dE_p}{d\sigma}=
\sum\limits_{|J|=p}[
\int\limits_{{\mathbb
R}^n}(D^{J}
P D^J\partial_\sigma Pd{\bf y} +
\int\limits_{{\mathbb
R}^n}D^J s
D^J\partial_\sigma s d{\bf y}+$$ $$ +\int\limits_{{\mathbb
R}^n}D^J
{\bf W} \Psi^{-1}(S,\sigma) D^J\partial_\sigma {\bf W} d{\bf y}
+$$
$$ +\int\limits_{{\mathbb
R}^n}(D^J {\bf W})^*
(\frac{\partial\Psi^{-1}(S,\sigma)}{\partial\sigma}- \frac{1}{\gamma}
\Psi^{-1}(S,\sigma)\partial_\sigma S) D^J {\bf W} d{\bf y}]= $$
$$ =I_1+I_2+ I_3 +I_4,$$ the integrals $I_k,\, k=1,2,3,4,$ are
numbered in the consecutive order.

We estimate every of the integrals in the standard way using the
H\" older and the Galiardo-Nirinberg inequalities (see
\cite{Serre},\cite{Grassin}) under assumption (2.3.11 -- 2.3.16) of
Theorem 1.  We denote by $c_i$ certain positive constants, which do
not depend on $\sigma.$

Let us begin from $I_1.$
$$I_1=
-\int\limits_{R_n}D^\alpha P D^\alpha(({\bf W},{\nabla}_{\bf y})P+
\frac{\gamma-1}{2}(P+\bar P) {\rm div}_{\bf y} {\bf W}+ ({\bf
W},{\nabla}_{\bf y})\bar P+PQ_1(t(\sigma)))d{\bf y}.$$
The integral from the terms of higher order
$\displaystyle D^\alpha P D^\alpha(({\bf
W},{\nabla}_{\bf y})P+ \frac{\gamma-1}{2}P {\rm div}_{\bf y} {\bf
 W})$ can be reduced by integration by parts to
$$ -\int\limits_{R_n}\left\{\frac{1}{2}|D^\alpha
P|^2 {\rm div}_{\bf y} {\bf W} + D^\alpha P\left(({\bf
W},{\nabla}_{\bf y}(D^\alpha P))- D^\alpha({\bf W},{\bf \nabla}_{\bf
y}P)\right)\right\} d{\bf y}-$$
$$ -\frac{\gamma-1}{2}
\int\limits_{R_n}D^\alpha P\left\{\left ({\nabla}_{\bf y}P,
D^{\alpha}{\bf W})+(
P {\rm div}_{\bf y}
 D^{\alpha}{\bf W}-D^{\alpha}(P {\rm div}_{\bf y}
 {\bf W})\right )\right\}
d{\bf y}
:=I_{11}.
$$
One can show applying the H\" older and the Galiardo-Nirinberg inequalities
that
$$ I_{11}\le
\tilde c_1\|D^p P\|_2(\|D^p P\|_2\|\nabla_ {\bf y}{\bf W}\|_\infty+
\|D^p {\bf W}\|_2\|\nabla_{\bf y} P\|_\infty)\le c_1
R^{1/2}F_m^{1/2}E_p $$ (the details in \cite{Serre},\cite{Grassin}.)
Further, according to  (2.3.5), (2.3.13) and (2.3.14)
$$\left|\int\limits_{R_n}D^\alpha P D^\alpha ({\bf W}{\nabla}_{\bf y}
\bar P)d{\bf y}\right|\le \tilde c_2
R^{1/2}\lambda^{-q}\exp(-\frac{\gamma-1}{2}{\rm tr} A)F_m \le c_2
R^{1/2}F_m,$$  $$\left|\int\limits_{R_n}D^\alpha P D^\alpha (\bar P)
{\rm div}_{\bf y} {\bf W})d{\bf y}\right
|\le c_3 R^{1/2}F_m,$$
$$\int\limits_{R_n}D^\alpha P
D^\alpha (PQ_1)d{\bf y}\ge 0.$$

The integral $I_2$ can be estimated analogously:
$$I_2\le c_4 R^{1/2}F_m^{1/2}E_p+ c_5 R^{1/2}F_m.$$
Further,
$$ I_3= \int\limits_{R_n} (D^\alpha{\bf W})^*\Psi^{-1} (S,\sigma)
\left \{({\bf W},{\nabla}_{\bf y}) D^\alpha{\bf W} -D^\alpha(({\bf
W}, {\nabla}_{\bf y}{\bf W})\right \} d{\bf y}+ $$ $$
+\frac{1}{2}\int\limits_{R_n} \left\{(D^\alpha{\bf W})^*\Psi^{-1}
(S,\sigma) D^\alpha{\bf W}) {\rm div} _{\bf y} {\bf W} - \frac{({\bf W
\nabla}_{\bf y}S)}{\gamma}(D^\alpha{\bf W})^* \Psi^{-1}(S,\sigma)
D^\alpha{\bf W}\right\} d{\bf y}+ $$ $$ +\frac{\gamma-1}{2}
\int\limits_{R_n}
\left\{(D^\alpha{\bf W})^*\Psi^{-1}(S,\sigma)
D^\alpha\left (P\Psi(S,\sigma)
{\bf
\nabla}_{\bf y}P + \Psi (S,\sigma)
{\bf \nabla}_{\bf y} (P\bar P) \right) \right\} d{\bf y}+
$$
$$
+\frac{\gamma-1}{4}
\int\limits_{R_n}
\left\{(D^\alpha{\bf W})^*\Psi^{-1}(S,\sigma)
D^\alpha\left (
(\Psi (S,\sigma)-\Psi (\bar S,\sigma)){\bf
\nabla}_{\bf y} (\bar P^2)
\right)
\right\}
d{\bf y}+ $$ $$
-\int\limits_{R_n}
(D^\alpha{\bf W})^*\Psi^{-1}(S,\sigma)
D^\alpha\left (Q_2{\bf W}-G\right )
d{\bf y}.
$$
First two integrals can be estimated from above by the values
$$c_6 R^{-1}\|D^p {\bf
W}\|_2^2\|\nabla {\bf W}\|_\infty+ $$$$c_7 R^{-1}\|D^p {\bf
W}\|_2^2\|{\bf W}\|_2 \le c_8 R^{1/2}F_m^{1/2}E_p,$$
the third, containing the
terms of the form $$D^p{\bf W}(\prod_{j=1}^k
(D^j S)^{\beta_j})D^{l-k}P D^{p+1-l}P,\,D^p{\bf W}(\prod_{j=1}^k (D^j
S)^{\beta_j})D^{l-k}\bar P D^{p+1-l}P,$$ or $$D^p{\bf
W}(\prod_{j=1}^k (D^j S)^{\beta_j})D^{l-k}P D^{p+1-l}\bar P,\,1\le
l\le p,\,1\le k \le l,\,\sum_j j\beta_j=k,$$ as follows from
the Galiardo-Nirinberg inequality and
(2.3.14), can be estimated by the sum of form $$R^{1/2}
E_p^{1/2}\sum_{i=9}^{11} c_i F_m^{1+\epsilon_i},\epsilon_i>0.$$
The fourth integral contains the term of form
$$\exp(\theta s) D^p{\bf W}(\prod_{j=1}^k (D^j
S)^{\beta_j})D^{l-k}\bar P D^{p+1-l}\bar P,\,1\le l\le p,\,1\le k \le
l,$$$$\sum_j j\beta_j=k,|\theta|<1,$$ and can be estimated in a
similar way, taking into account that the value of $\|s\|_\infty$ is
bounded.

The last integral taking into account
condition (2.3.16) can be estimated by the quantity $c_9 F_m.$
Remark here that
$$G=\lambda^{-1}(t)A_1(t)\nabla_{\bf V}{\bf f}_1
(t,{\bf x},\bar\Pi+\theta_1 \pi,\bar {\bf V}+\theta_2 {\bf v},
\bar S+\theta_3 s)A_1^{-1}(t){\bf W}+$$$$\lambda^{q-2}(t)
\nabla_\Pi{\bf f}_1
(t,{\bf x},\bar\Pi+\theta_1 \pi,\bar {\bf V}+
\theta_2 {\bf v},
\bar S+\theta_3 s)P+$$$$\lambda^{-2}(t)
\nabla_\Pi{\bf f}_1
(t,{\bf x},\bar\Pi+\theta_1 \pi,\bar {\bf V}+\theta_2 {\bf v},
\bar S+\theta_3 s)s
,$$where $\theta_i\in(0,1),\,i=1,2,3,$
and
$$
\int\limits_{R_n}
(D^\alpha{\bf W})^*\Psi^{-1}(S,\sigma)
D^\alpha(U_\phi(t){\bf W})
d{\bf y}=0,
$$ since
$$({\bf w},A_2^{-1}(t)U_\phi(t){\bf w})=0$$
for any real-valued vector ${\bf w}.$

 Further we estimate $I_4.$
Remark first of all, that according to (2.3.13)
$R'(\sigma)/R(\sigma)\ge 0.$ Taking into account (2.3.15)
we obtain
$$\frac{\partial\Psi^{-1} (S,\sigma) }{\partial
\sigma}=-\frac{\Psi^{-2}(S,\sigma) \Psi(S,\sigma) R'}{R}+
O(\Psi^{-1}(S,\sigma) )\le c_{12}\Psi^{-1}(S,\sigma).$$ Thus,
$$I_4\le \tilde c_{13}(1+\|{\bf W}{\bf \nabla} _{\bf y} s\|_\infty+
\|{\bf W}{\bf \nabla}_{\bf y} \bar S\|_\infty)E_p\le c_{13}(1+
R^{1/2} F_m+ R^{1/2} F_m^{1/2} ) E_p.$$

So we get $$ I_{1}\le c_1(
E_p+R^{1/2}F_m^{1/2}E_p+ R^{1/2}F_m),$$ $$I_2\le
c_2(R^{1/2}F_m^{1/2}E_p+ R^{1/2}F_m),$$ $$ I_3\le c_3(
F_m+R^{1/2}F_m^{1/2}E_p+ R^{1/2} F_m+R^{1/2}E
_p^{1/2}\sum_{i=1}^4 F_m^{1+\gamma_i}),\gamma_i>0,$$
$$I_4\le c_5(1+ R^{1/2} F_m+ R^{1/2} F_m^{1/2} ) E_p.$$
Using all the
estimates we get  $$F_m'\le c_6(F_m+R^{1/2}(F_m+
\sum_{i=1}^6 F_m^{3/2+\tilde\gamma_i}),\tilde\gamma_i\ge 0.$$

Set $\Lambda_m(\sigma)=e^{-{c_6}\sigma}F_m.$
Then $$\Lambda_m'\le c_7 (\Lambda_m +
\sum_{i=1}^6
(\Lambda_m)^{3/2+\tilde\gamma_i})R^{1/2},\tilde\gamma_i\ge
0,\eqno(2.3.17)$$ where the constant $c_7$ depends on
$c_6,\,\sigma_\infty$.
Let
$$\Theta(\tilde g):=\int\limits_{\delta}^{\tilde g}
\frac{dg}{g+\sum_{i=1}^6
g^{3/2+\tilde\gamma_i}},$$
$\delta>0.$ The integral
diverges in the zero, such that $\Theta(0)=-\infty.$
Integrating inequality (2.3.17) over
$\sigma$ we obtain $$ \Theta(\Lambda_m(\sigma))\le
\Theta(\Lambda_m(0))+ c_8\int\limits_{0}^\sigma
R^{1/2}(\tilde\sigma) d\tilde\sigma,$$
moreover, as follows from (2.3.12),
the integral in the right-hand side of the last inequality
converges as
$\sigma\to\sigma_\infty$ to the constant $C$, depending only on
the initial data.
Choosing $\Lambda_m(0)$ (that is the
$H_m -$ norm of the initial data) sufficiently small,
one can make the value $\Theta(\Lambda_m(0))+C$
later then $\Theta(+\infty),$ it signifies that $\Lambda_m(\sigma)$
and $F_m(\sigma)$ are bounded from above for all
$\sigma\in[0,\sigma_\infty) $ and $\sigma_*=\sigma_\infty$.  So,
${\bf u}\in \cap _{j=0}^m
C^j([0,\infty);H^{m-j}({\mathbb R}^n))$ and
Theorem 2.1 is proved.

{\sc Remark 2.3.2}
The solutions with linear profile of velocity,
satisfying the  conditions of Theorem 2.1, exist.  In the case
investigated in
\cite{Serre}, ${\bf f}=0,\,\gamma\le
1+\frac{2}{n},\,$
$A(t)=(E+tA(0))^{-1}A(0), \,$
$ {\rm Sp} A(0)\cap\mathbb R_-=\emptyset,\,$${\bf b}(t)=0,$
$A_1(t)=A(t),$
$\bar\rho=0, \,\bar S=const.$
Here
$\lambda=(1+t)^{-2},\,$$q=\frac{n(\gamma-1)}{4}.$

%\subsection{Interior solution with linear profile of velocity.}
\medskip
\subsubsection{Corollary and examples}

The following Corollary helps us to prove that among the solutions
to the Euler system constructed in section 2.1 and 2.2 there are
interior ones.
\theoremstyle{corollary 2.}
\newtheorem{corollary 2.}{Corollary 2.}
\begin{corollary 2.}
Let ${\bf f}=L(t){\bf V},$
where $L(t)=-\mu I+U_1(t),$
$\mu  $ is a nonnegative constant, $U_1(t)
$ is a one-parameter family of skew-symmetric matrices with the
smooth coefficients.  Assume that the system (2.3.1--2.3.3) has a
global in time classical solution $(\bar \rho, A(t){\bf r}, \bar S)$,
satisfying the condition a) of Theorem 1 and $A(t) \sim \alpha(t)
I+\beta(t)U_2,$ as $t\to \infty,$ where $\alpha(t) $ and $\beta(t)$
are the real-valued functions, $U_2$ is a skew-symmetric matrix with
constant coefficients.  Moreover, suppose that $\alpha(t)=
\frac{\delta}{t},$ with $\,\delta=const$, such that
$\delta>0$ for $\mu>0, $ and $\delta>\frac{1}{2}$ for $\mu=0$. Then
this solution is interior.  \end{corollary 2.}

{\it Proof.} Here ${\bf b}(t)=0.$
The matrix $A_1(t)=  \exp (-\int_{t_0}^t A(\tau)d\tau)$
satisfying the condition b) of Theorem 2.1
exists, $\det A_1(t)=\exp(-n \int_{t_0}^t
\alpha(\tau)d\tau),$ as the eigenvalues of a skew-symmetric
real-valued matrix have the real part equal to zero.  The first and
second conditions in (2.3.16) hold, since ${\bf f}$ is independent on
$\rho$ and $S$, $Q_3$ depends only on $t$(see Remark 2.3.1).  We
choose as $U_\phi(t)$ the  matrix
$A_2(t)A_1(t)(\beta(t)U_2-U_1(t))A_1^{-1}(t)$
(such that $A_2^{-1}(t)U_\phi(t)$ is skew-symmetric)
to obtain
$Q_3=\lambda^{-1}(t)((\ln \lambda(t))'+2\alpha(t)+\mu A_2(t))I.$
If $\lambda(t)=
\exp(-\varepsilon\int_{t_0}^t \alpha(\tau)d\tau), \varepsilon=const>0,
\varepsilon\delta >1,$ then both integrals in (2.3.11), (2.3.12)
converge.  To guarantee the nonnegativity of $R'(t), Q_1(t), Q_3(t)$
we must satisfy the following inequalities for $\varepsilon$ and $q$:
$q<1,\,$ $\varepsilon \ge \frac{1}{1-q}, \,
$$\varepsilon\le\frac{(\gamma-1)n}{2q},\, \varepsilon\le 2$
(the last inequality arises only in the case $\mu=0$), the
conditions (2.3.12) hold for $\varepsilon\le\frac{(\gamma-1)n}{2q}.$
For $\mu>0$, it implies $\frac{1}{1-q}\le \varepsilon \le
\frac{(\gamma-1)n}{2q},$ $\frac{1}{\delta}< a \le
\frac{(\gamma-1)n}{2q},\,$ $q<\bar q :=\min\{\frac{\delta n
(\gamma-1)}{2},\,\frac{n(\gamma-1)}{2+n(\gamma-1)}\},$ and we choose
$q\in (0,\bar q).$ If $\mu=0,$ then additionally we have
inequalities $\frac{1}{\delta}< \varepsilon \le 2,
\,\frac{1}{1-q}\le \varepsilon \le 2.$ Since $\delta<\frac{1}{2},$
then we can choose $\varepsilon$ to satisfy the first inequality,
for the second we have $q\le \frac{1}{2}.$ In this case we choose
$q\in (0,\min\{{\bar q},\frac{1}{2}\}).$ Thus, all conditions of
Theorem 2.1 are satisfied and the Corollary is proved.

\medskip
{\sc Examples.}
Thus, as follows from  Corollary 2.1, in the cases
$ \mu=l=0 $ and $ \mu>0 $ in the Sections 2.1.1 ($n=2$) and 2.2
($n=3$)
we
have constructed the velocity field for the interior solution.
Really, the right-hand side has the form indicated in the
Corollary statement, the velocity has linear profile ${\bf
V}=A(t){\bf r},$ with $A(t)$ of special form, that is $A(t)=\alpha(t)I
+\beta(t)U_2$ with the skew-symmetric $U_2.$
For $\mu=l=0,$ we have $\alpha(t)\sim t^{-1},
t\to\infty,\,\delta=1>\frac{1}{2}.$

For $\mu>0, $ we get
$\alpha(t)\sim\frac{1}{2\gamma}t^{-1},\, t\to\infty,
\delta=\frac{1}{2\gamma}>0$ ($n=2$), and
$\delta=\frac{1}{3\gamma-1}>0$ ($n=3$)(remark that the
denotation $\delta$ has here other sense than in the subsection 2.2.)

The solution constructed in the Section 2.1.2 is also interior
for $\mu=l=0$. As follows from Proposition 2.1, in the case
$A(t)\sim \alpha(t) I,$ where
$\alpha(t)= \delta t^{-1}, t\to\infty,\,\delta>\frac{1}{2}.$

Moreover, as the numerical analysis suggests,
all solutions with linear profile of velocity (of general form) are
interior if $ \mu>0,\,l=0,$ either for ${\bf b}(t)={\bf 0}$ or ${\bf
b}(t)\ne {\bf 0}$ (it seems that in the last case ${\bf b}(t)\sim
{\bf B}t^{-1}, $ with the constant vector ${\bf B}$).  However, the
last assertion can be considered only as a hypothesis.
The case
$\mu>0,
\,l\ne 0$ is more complicated, because it seems that it accepts
stable equilibriums except of the origin.

{\sc Remark 2.3.3}
We can find the entropy function for the solution constructed in
Sections 2.1 and 2.2. For example, in the simplest case (Section
2.2.1)
$$ S(t,|r|, \phi)= S_0(|r|\exp(-\int\limits_0^t \alpha(\tau)
d\tau), \phi+\int\limits_0^t \beta(\tau) d\tau).  $$
For example, for the initial data (2.1.12)
$S_0=const+(a(\gamma-1)+\gamma)\ln(1+|{\bf r}|^2).$

In
spite of the components of the pressure and the density must vanish
as $|{\bf x}|\to\infty,$ the entropy may even increase (remember
that the conditions to Theorem 2.1 require  only the boundedness of
the entropy gradient).

{\sc Remark 2.3.4} In the case $G(0)\ne 0$ the
density cannot be compactly supported
(in contrast with $G(0)=0$).  Really, since
$p=\pi^\frac{2\gamma}{\gamma-1},\, \rho=
\pi^\frac{2}{\gamma-1}e^{-\frac{S}{\gamma}},$ then due to the
compatibility conditions
$\pi\nabla\pi\sim const\cdot e^{-\frac{S}{\gamma}}, \,|{\bf x}|\to
c-0,$ where $c $ is a point of the support of $\pi.$
Therefore, for $C^1 -$ smooth $\pi$ it occurs that $S\to
+\infty, \,|{\bf x}|\to c-0,$ and we cannot choose
any smooth initial data.  It is interesting that
if one requires only $C^0 -$ smoothness of
$\pi\,(\pi\sim const\cdot (c-|{\bf x}|)^{1/2},\,|{\bf x}|\to c-0, \,
\pi=0,\,|{\bf x}|\ge c>-0),$ the condition may be fulfilled.
Moreover, $\rho$ and $p$ will be of the $C^1$- class of smoothness,
however, neither the theorem on the local in time existence of the
smooth solution, no Theorem 2.1 can be applied.

\section{Acceptable velocities:\
generalization of velocity with
linear profile}

\subsection
{The auxiliary system of transport equations}

Consider the following system of equations associated with (E1--E3):
$$ \partial_t
{\bf V}+({\bf V},{\bf \nabla})\,{\bf V}= 0.\eqno(3.1.1)$$
Here
$\bf V$ is a vector field from the tangent bundle of $\Sigma.$

We look for the smooth solution to (3.1.1)
of the separated form
$${\bf V}=A(t){\bf\Lambda(x)}\eqno (3.1.2)$$
with a square $(n\times n)$ matrix $A(t).$

The following possibilities can be realized:

{\bf I.}
$$\nabla_i\Lambda^j=const \cdot\delta_i^j,\eqno(3.1.3)(A1)$$
with the Kronekker symbol
$\delta_i^j,$
$A'(t)=-A^2(t).$

If $Sp(A(0))\cap{\mathbb R}_-=\emptyset,$ then $A(t)$ is bounded all
over the time $t\ge\infty.$

It is easy to see that for the euclidean metrics there are the
unique possibility $\Lambda=a{\bf r}+b, $ where $\bf r$ is the
radius-vector of point, $a, b$ are constants, so
we come back to the velocity with the linear profile.

{\bf II.}
$$\Lambda^i=\Lambda^j\nabla_j \Lambda^i, \,
i=1,...,n,\eqno(3.1.4)(A2)$$
and
$A(t)=a(t)I,$  with the identity matrix $I$,
$a(t)$ is a function such that
$a'(t)=-a^2(t).$

In the second case we essentially restrict the set of matrices,
but the set of corresponding vector fields $\bf\Lambda$ is more rich
then in the case (A1).

{\sc Remark 3.1.1.} The equations (A1) may be incompatible (f.e. for
the 2D sphere with natural coordinates).

{\sc Remark 3.1.2.} If $\bf\Lambda$ satisfies (A1), then it satisfies
(A2).

\subsubsection{The vector field satisfying condition (A2)
with a constant divergency.}

In our further considerations the vector fields $\bf\Lambda$
satisfying (A2) having a constant divergency, play an important
role, therefore we study them in detail.
\theoremstyle{theorem 3.}
\newtheorem{theorem 3.}{Theorem 3.}
\begin{theorem 3.}
If the nontrivial potential vector field satisfying
(A2) has a constant divergency
$D=D_0$, then $D_0>0.$
\end{theorem 3.}

{\it Proof.}
Really, $$({\bf\nabla
\Lambda})=\sum\limits_{i,j}\nabla_i(\Lambda^j\nabla_j\Lambda^i)=
2\sum\limits_{i\ne
j}\nabla_i\Lambda^j
\nabla_i\Lambda^j
+\sum\limits_{i}(\nabla_i\Lambda^i)^2+
\sum\limits_{j}\Lambda^j\nabla_j({\bf\nabla\Lambda}).\eqno(3.1.5)
$$
(in the formula we do not sum over indices $i$ in
$(\nabla_i\Lambda^i)$).

If the field $\bf\Lambda$ is potential (that is there exist
a function $\Phi$ such that ${\bf\Lambda}={\bf\nabla}\Phi$),
then $\nabla_i\Lambda^j=\nabla_j\Lambda^i$. If we suppose that
the divergency is constant, then $$D_0=
2\sum\limits_{i\ne
j}(\nabla_i\Lambda^j)^2+\sum\limits_{i}(\nabla_i\Lambda^i)^2.  $$
Thus, $D_0$ is a sum of squares and it is therefore nonnegative.
However, if $ D=\bf\nabla
\Lambda=0,$ then $\nabla_i \Lambda^j=0$ under all combinations of
indices, therefore $\bf\Lambda$ is constant and,
as follows from (A2), equal to zero, that is trivial in spite of
the assumption.

Theorem 3.1 has the evident corollary.

\begin{cor 3.}
There is no nontrivial potential divergency free vector field
satisfying (A2).
\end{cor 3.}

{\sc Remark 3.1.3.} As we shall see below, the possibility to choose
the divergency free vector field should involve the possibility to
construct the solutions to the Euler equations with the properties of
localization of mass in a point. But the situation is
restricted for the systems, containing pressure.

The case of dimension greater then one is more interesting for us,
as for the dimension equal to one the situation is the
following:  if nontrivial vector field satisfies condition (A2),
then its divergency is equal to 1 automatically.

Further we assume that the field
$\bf\Lambda$ has at least two times differentiable components.

Let us remark that a smooth vector field with constant nonzero
divergency cannot exist on the compact manifold without boundary.
Indeed, it contradicts to the fact that the integral from the
divergency taken over all manifold has to be zero.

\begin{theorem 3.}
There exists no nontrivial divergency free vector field given on the
manifold of dimension $n$ satisfying condition (A2).
\end{theorem 3.}

Denote $$J_m=\sum\limits_{i_k\in K_m,\,j_k\ne i_k,
\,k=1,...,m} (\nabla_{i_1}\Lambda^{i_1}...\nabla_{i_m}\Lambda^{i_m}-
\nabla_{j_1}\Lambda^{i_1}...\nabla_{j_m}\Lambda^{i_m})
,$$
where the summation is taken over all subsets
$K_m={i_1,...,i_m}, \,1\le m\le n,$ from the set
${1,...,n}$ and the set ${j_1,...,j_m}$ is a permutation
of elements of $K_m.$

We need the following lemma
\newtheorem{lemma 3.}{Lemma 3.}
\begin{lemma 3.}
Let the field $\bf\Lambda$ satisfy (A2) and
have a constant divergency $D.$ Then
the following representation holds:
$$D=D^m+const\cdot J_m+ {\mathcal
M}(D,J_2,...,J_{m-1}),\eqno(3.1.6)$$ where ${\mathcal
M}(D,J_2,...,J_{m-1}) $ is a polynomial of $k$-th degree $0<k<m$,
$\,m=2,...,n.$ \end{lemma 3.}

{\it Proof of Lemma 3.1.}  From (A2) we get
$$\Lambda^{i_1}=\Lambda^{i_2}\nabla_{i_2}\Lambda^{i_1}=
\Lambda^{i_3}\nabla_{i_3}\Lambda^{i_2}\nabla_{i_2}\Lambda^{i_1}=
\Lambda^{i_m}\nabla_{i_m}\Lambda^{i_{m-1}}...\nabla_{i_2}\Lambda^{i_1},
$$
where $i_k=1,...,n,\,k=1,...,n,\,m\le n.$
Therefore
$$D=\sum\limits_{i_1=1}^{n}\nabla_{i_1}\Lambda^{i_1}
=
\nabla_{i_1}\Lambda^{i_m}
\nabla_{i_m}\Lambda^{i_{m-1}}...\nabla_{i_2}\Lambda^{i_1}+
\Lambda^{i_m}\nabla_{i_1}
(\nabla_{i_m}\Lambda^{i_{m-1}}...\nabla_{i_2}\Lambda^{i_1})+
$$$$
\nabla_{i_1}\Lambda^{i_m}
\nabla_{i_m}\Lambda^{i_{m-1}}...\nabla_{i_2}\Lambda^{i_1}+
+\Lambda^{i_m}
\nabla_{i_1}(\nabla_{i_m}\Lambda^{i_1})=
 $$$$\nabla_{i_1}\Lambda^{i_m}
\nabla_{i_m}\Lambda^{i_{m-1}}...\nabla_{i_2}\Lambda^{i_1}+
\Lambda^{i_m} \nabla_{i_m}D.$$ If $D=const,$ then $D=D^m+const\cdot
J_m+ {\mathcal M}(D,J_2,...,J_{m-1}),$ the last polynomial is
homogeneous of $m$-th degree with respect to
$\nabla_i\Lambda^j,\,i,j=1,...,n.$

\medskip
{\sc Example.} The rather thorough
(for
$m>2$)
calculation shows that the following formulas are true:
$$D=D^2-2J_2,\eqno(3.1.7)$$ $$D=D^3+3J_3-3DJ_2, \eqno(3.1.8)$$
$$D=D^4+4J_4-4D^2J_2+DJ_3+2J_2^2.\eqno(3.1.9)$$

\begin{lemma 3.}
For the field $\bf\Lambda$ satisfying (A2)
the following representation takes place
$$1-D+J_2-J_3+J_4-...+J_n=0,\eqno(3.1.10)$$
\end{lemma 3.}

{\it Proof of Lemma 3.2.} Consider the system,
defined by the condition (A2) as a linear homogeneous system
with respect to the variables $\Lambda^i$.
The condition (3.1.10) is namely the condition of existence
of its nontrivial solution, that is
the corresponding determinant of linear system is equal to zero.
So the proof is over.

{\it Proof of Theorem 3.2.} Suppose that $D=0.$ Then from the
representation (3.1.6) we obtain consecutively that
$J_m=0,\,m=2,...,n.$  This fact contradicts to (3.1.10).

\begin{theorem 3.}
Suppose that on the manifold $\Sigma$ of dimension $n$ there exists
the vector field $\bar{\bf\Lambda}$ satisfying condition (A1).
If a nontrivial vector field, given on $\Sigma$ and  satisfying
condition (A2) has a constant divergency $D$, then $D$ can take only
integer values from 1 to n.  \end{theorem 3.}

{\it Proof.} Equations (3.1.6) for $m=2,...,n$ and (3.1.10)
form the system of  $n$ algebraic equations for
$n$ unknown variables $(D,J_2,...,J_n)$.  We express consecutively
$J_2,...,J_n$ from equations (3.1.6) through the divergency
$D$ and substitute the result in (3.1.10). In such way we obtain the
equation of $n$-th degree  for $D$. Thus, $D$ cannot take more then
$n$ different real values.  But the field $\bar{\bf\Lambda},$ as well
as its projections to the subspaces of dimension $1,...,n-1$, satisfy
condition (A2), moreover, their divergency takes exactly $n$ values,
namely, $1,2,...,n.$ Thus, $D$ cannot take another values and the
Theorem 3.3 is proved.

{\sc Remark 3.1.4.} For the euclidean space
$\bar{\bf\Lambda}=\bf r,$ where $\bf r$ is the radius-vector of
point.

{\sc Remark 3.1.5.} To all appearances, the requirement of existence
of the field satisfying (A1), is unnecessary, that is in the case
where the divergency of the field $\bf \Lambda$ with the property
(A2) is constant, it is equal to a natural number, later or equal to
the dimension of space. At least the fact holds in the spaces of
dimensions 2 and 3, as the following proposition shows.

\begin{utv 3.}
If a nontrivial vector-field given on the manifold
$\Sigma$ of dimension $n=2$ or $n=3$ and satisfying
condition (A2) has the constant divergency $D$, then $D$
can take only
natural values from 1 to n.
\end{utv 3.}

{\it Proof.} a) If $n=2,$ then from (3.1.6), (3.1.10) we have that
$D=D^2-2J_2,\,1-D+J_2=0,$ or $D^2-3D+2=0$. The last equation
has two roots $D=1,\,D=2.$

b) Analogically for $n=3$ we have the system of equations
$D=D^2-2J_2,\,D=D^3+3J_3-3DJ_2,\,1-D+J_2-J_3=0.$
Therefore $D^3-6D+11D-6=0,$ and the roots are
$D=1,\, D=2,\, D=3.$

\subsubsection{A generalization of condition
(A2)}

The solution to the equation (A2) is not unique. Suppose that we
succeed in finding of finite or infinite number of the solutions
such that ${\bf\Lambda_k,\, \bf k\in K,\, K}\subseteq \mathcal N,$
having the additional property
$$ \Lambda^j_{\bf m}\nabla_j \Lambda^i_{\bf l}+\Lambda^j_{\bf
l}\nabla_j \Lambda^i_{\bf m}= \sum\limits_{\bf k} \beta_{\bf
k}{\bf\Lambda_k}, \, {\bf m}\ne {\bf l}, \, \beta_{\bf k}=const.,{\bf
k,m,l}\in {\bf K}.\eqno(3.1.11)$$ Then there exists a solution
of the form $${\bf u}=\sum\limits_{\bf k} a_{\bf
k}(t){\bf\Lambda_k}.$$ The functions $a_{\bf k}(t)$
satisfy the system of equations
$$a'_{\bf k}(t)=-a^2_{\bf k}(t)-\beta_{\bf k}a_{\bf
k}(t)\sum\limits_{{\bf l}\ne {\bf k}} a_{\bf l}(t),\,\bf k\in
K.\eqno(3.1.12)$$

\section{Two-dimensional manifold}

The section is devoted to the important for applications case of
a two dimensional manifold.

\subsection{Transport equation on the manifold covered by one
chart}

Let $(x^1, x^2)$ be the coordinates on the manifold,
and $(\Lambda^1, \Lambda^2)$ be the components of the solution to
equation (3.1.4).  Thus, the equations (A2) has the form $$
\Lambda^1=\Lambda^1\frac{\partial \Lambda^1} {\partial x^1}+\Lambda^2
\frac{\partial \Lambda^1} {\partial
x^2}+\Gamma^1_{11}(\Lambda^1)^2+2\Gamma^1_{12}
\Lambda^1\Lambda^2+\Gamma^1_{22}(\Lambda^2)^2,\eqno(4.1.1)$$
$$ \Lambda^2=\Lambda^1\frac{\partial
\Lambda^2}
{\partial x^1}+\Lambda^2 \frac{\partial
\Lambda^2}
{\partial x^2}+\Gamma^2_{11}(\Lambda^1)^2+2\Gamma^1_{12}
\Lambda^1\Lambda^2+\Gamma^2_{22}(\Lambda^2)^2\eqno(4.1.2)$$

If we denote $Z=\frac{\Lambda^2}{\Lambda^1},$
then as a corollary of (4.1.1--4.1.2) we obtain
the quasilinear equation
$$\frac{\partial
Z}{\partial x^1}+Z\frac{\partial Z}{\partial x^2} = \Gamma^1_{22}
(Z)^3+(2\Gamma^1_{12}-\Gamma^2_{22})(Z)^2+
(\Gamma^1_{11}-2\Gamma^2_{12})Z-\Gamma^2_{11},\eqno(4.1.3)$$
which can be solved by the characteristics method:
$$\frac{d Z}{dx^1}= \Gamma^1_{22}
(Z)^3+(2\Gamma^1_{12}-\Gamma^2_{22})(Z)^2+
(\Gamma^1_{11}-2\Gamma^2_{12})Z-\Gamma^2_{11},$$
$$\frac{d x^2}{d
x^1}=Z.$$

If  $Z(x^1,x^2)$ is found, then from (4.1.1) we obtain the
linear equation $$ \frac{\partial
\Lambda^1} {\partial x^1}+Z\frac{\partial \Lambda^1} {\partial
x^2}+(\Gamma^1_{11}+2\Gamma^1_{12}Z+
\Gamma^1_{22}(Z)^2)\Lambda^1=1,\eqno(4.1.4)$$
for the function $\Lambda^1.$ As
$\Lambda_2=Z\Lambda_1,$ then the solution components are found.

If $\lambda_i=\sqrt{g_{ii}}\Lambda^i$ are the physical components
of the field $\bf \Lambda,$
then $$ \frac{\partial \lambda_1} {\partial
x^1}+\frac{\sqrt{g_{11}}}{\sqrt{g_{22}}}\frac{\lambda_2}{\lambda_1}
\frac{\partial \lambda_1} {\partial
x^2}+\frac{\sqrt{g_{11}}}{\sqrt{g_{22}}}\Gamma^1_{12}\lambda_2+
\frac{{g_{11}}}{{g_{22}}}\Gamma^1_{22}\frac
{(\lambda_2)^2}{\lambda_1}=\sqrt{g_{11}},$$

$$ \frac{\partial
\lambda_2}
{\partial
x^1}+\frac{\sqrt{g_{11}}}{\sqrt{g_{22}}}\frac{\lambda_2}{\lambda_1}
\frac{\partial \lambda_2} {\partial
x^2}+\frac{\sqrt{g_{22}}}{\sqrt{g_{11}}}\Gamma^2_{11}\lambda_1+
\Gamma^2_{21}
\lambda_2=\sqrt{g_{11}}\frac
{\lambda_2}{\lambda_1}.$$

If $z=\frac
{\lambda_2}{\lambda_1}$,
then the function satisfies the equation
$$ \frac{\partial z} {\partial
x^1}+\frac{\sqrt{g_{11}}}{\sqrt{g_{22}}}z
\frac{\partial z} {\partial
x^2}+\frac{\sqrt{g_{22}}}{\sqrt{g_{11}}}\Gamma^2_{11}+
\Gamma^2_{21}
z- \frac{\sqrt{g_{11}}}{\sqrt{g_{22}}}\Gamma^1_{12}z^2-
\frac{g_{11}}{g_{22}}\Gamma^1_{22}z^3=0,\eqno(4.1.5)$$
the components $\lambda_1$ and $\lambda_2$ can be found from
the equation
$$ \frac{\partial \lambda_1}
{\partial x^1}+\frac{\sqrt{g_{11}}}{\sqrt{g_{22}}}z\frac{\partial
\lambda_1} {\partial
x^2}+\frac{\sqrt{g_{11}}}{\sqrt{g_{22}}}\Gamma^1_{12}{\lambda_1}{z}+
\frac{{g_{11}}}{{g_{22}}}\Gamma^1_{22}z^2
\lambda_1=\sqrt{g_{11}},$$
$$ \frac{\partial
\lambda_2}
{\partial
x^1}+\frac{\sqrt{g_{11}}}{\sqrt{g_{22}}}z\frac{\partial \lambda_2} {\partial
x^2}+\frac{\sqrt{g_{22}}}{\sqrt{g_{11}}}\Gamma^2_{11}\frac{\lambda_2}{z}+
\Gamma^2_{21}
\lambda_2=\sqrt{g_{11}}z.\eqno(4.1.6)$$
Here
$\lambda_1=\displaystyle\frac{\lambda_2}{z}.$

\medskip
{\sc Remark 4.1.1.} As we shall see, the important role plays a
possibility to find a potential vector field satisfying (A2).
In the terms of the  potential $\Phi$ (such that ${\bf \Lambda}={\bf
\nabla}\Phi$)  condition (A2) takes the form
$${\bf \nabla}\Phi=({\bf \nabla}\Phi,{\bf \nabla}){\bf \nabla}\Phi$$
or
$$\nabla^1\Phi=\nabla^1\Phi\nabla_1^1\Phi+\nabla^2\Phi\nabla^1_2\Phi,$$
$$\nabla^2\Phi=\nabla^1\Phi\nabla^2_1\Phi+\nabla^2\Phi\nabla^2_2\Phi.$$

If we consider the system as linear homogeneous with respect to
$\nabla^1\Phi$ and $\nabla^2 \Phi$ we obtain the necessary condition
for the existence of the nontrivial potential $\Phi:$
$${\rm Hess} \Phi-\Delta \Phi+1=0,$$
where
${\rm Hess}
\Phi=\nabla^1_1\Phi\nabla_2^2\Phi-\nabla_1^2\Phi\nabla_2^1\Phi,\quad
\Delta \Phi=\nabla^1_1\Phi+\nabla_2^2\Phi.$

 \medskip {\sc Remark 4.1.2.} For the equation of the form
$$ \Xi (x^1,x^2,{\bf \Lambda})\Lambda^i=\Lambda^j\nabla_j \Lambda^i,
\quad i=1,2,$$ whose particular case is (4.1.2), the equation (4.1.3)
is true as well.  But instead of (4.1.4) we obtain the equation $$
\frac{\partial \Lambda^1} {\partial x^1}+Z\frac{\partial \Lambda^1}
{\partial x^2}+(\Gamma^1_{11}+2\Gamma^1_{12}Z+
\Gamma^1_{22}(Z)^2)\Lambda^1=\Xi(x^1,x^2,|{\bf
\Lambda}|).\eqno(4.1.7)$$

In the case of aerodynamic friction, where
${\bf F}({\bf V})=-\mu_1{\bf V}|{\bf V}|,\,\mu_1=const>0,$
we can obtain
$$\Xi(x^1,x^2,|{\bf \Lambda}|)=1-\mu_1\,{\rm sign}\,
a(t)|\Lambda^1|\sqrt{1+(Z)^2}.$$

\medskip
{\sc Remark 4.1.3.} If the dimension $n$ is greater than 2, by
means of introduction of a new variable
$Z_{k-1}=\frac{\Lambda^k}{\Lambda^1},\,k=2,...,n,$
acting analogously,
one can deduce the equation (A2) to the
system of $n-1$ equations for $Z_i,\,i=1,...,n-1.$
In the case of euclidean metrics it will be
$n-1 -$ dimensional system of transport equations as well,
therefore the procedure of the dimension reduction  we can
apply several times up to reduction to one equation.

\subsubsection{Nonlinear transport equation with
the Coriolis force on a two dimensional surface}

We mean the equation

$${\bf u}_t+({\bf u,\nabla}){\bf u}=l{\bf u}_{\bot},$$
or
$$u^i_t+u^j\nabla_j u^i=le_{\cdot j}^i u^j,\eqno(4.1.8)$$
where $l=l(x^1,x^2)$ is a Coriolis parameter,
$u_{\bot}^i=e_{\cdot j}^i u^j, \,e_{ij}$ is the skew-symmetric
Levi-Civita tensor.

We seek for a solution of the equation in the form
$${\bf u}= a(t){\bf\Lambda}+{\bf
\Xi},\eqno(4.1.9)$$ where $\bf\Lambda$ and $\bf\Xi$ are vectors,
which do non depend on time.

In the case  $a'(t)=-a^2(t)$ as before, the vector
$\bf\Lambda$ satisfies condition (A2), the vector $\bf\Xi$
is a stationary solution to (4.1.8), that is the solution to the
equation
$$le_{\cdot j}^i\Xi^j=\Xi^j\nabla_j \Xi^i, \,
i=1,...,n.\eqno(4.1.10)$$

Moreover, the vectors $\bf\Lambda$ and
$\bf\Xi$ can be compatible in some sense, that is have to satisfy the
equation
$$le_{\cdot
 j}^i\Lambda^j=\Xi^j\nabla_j \Lambda^i+ \Lambda^j\nabla_j \Xi^i, \,
i=1,...,n.\eqno(4.1.11)$$

Thus the solution of form (4.1.9)
is a sum of a nonstationary solution to the non rotational
equation (4.1.1) and
a stationary solution to (4.1.8).
Of course, it is not evident before that the condition
(4.1.11) takes place.

\subsection {The case of the plane}

For the euclidean metrics all Christoffel symbols are equal to zero,
the geometrical components of the vector $\bf \Lambda$
are equal to
physical ones.  Therefore the equation (4.1.3) coincides with the
transport equation (here $(x_1, x_2) =(x^1,x^2) $ are the charthesian
coordinates):  $$ \frac{\partial z}{\partial x_1}+z\frac{\partial
z}{\partial x_2} =0.\eqno(4.1.12)$$ As well known (f.e.,
\cite{Whitham}), the solution to equation (4.1.12) can be given
implicitly in the form $z=F(x_2 - z x_1),$ where $z(0, x_2)= F(x_2)$
is the Cauchy datum on the noncharacteristic curve $x_1=0.$ The
functions $\Lambda_1(x_1,x_2), \Lambda_2(x_1,x_2)$ can be found from
equations $$ \frac{\partial \Lambda_1}{\partial x_1}+z\frac{\partial
\Lambda_1}{\partial x_2} =1,$$ $$\frac{\partial \Lambda_2}{\partial
x_1}+z\frac{\partial \Lambda_2}{\partial x_2} =z.$$
Note that the
presence of singularities of a solution to (4.1.12) does not signify
the presence of singularities of $\Lambda_1, \Lambda_2.$

The solution with linear profile of velocity
of the form ${\bf V}=a(t)I{\bf
r}+ B(t),$ where ${\bf r}$ is a radius-vector of point,
corresponds to $z=\frac{a x_2+b_2}{a x_1+b_1}.$

It is not difficult to see that
$z=K=const$ corresponds to the velocity field
$$\Lambda_1=x_1+\phi(x_2-K
x_1),\,\Lambda_2=Kx_1+K\phi(x_2-Kx_1),\eqno(4.1.13)$$
where $\phi$
is an arbitrary differentiable function.
The divergency of the velocity field is equal to 1.

Remark that it can be potential only if it coincides (up to the
rotation with respect to the  origin)
with
the linear field, depending only on one coordinate.

Thus, we have found the method of constructing  solutions with
separated variables of the form
${\bf u}=a(t)\bf\Lambda$ of the two-dimensional transport
equation on the plane.

Stress that the solution cannot be divergency free.

Let us try to find the condition of the existence of a solution
of form $${\bf
u}=a_1(t){\bf\Lambda_{\bf 1}}+ a_2(t){\bf\Lambda_{\bf 2}},$$ where
${\bf\Lambda_{\bf 1}}$ is given by formula (4.1.13), and
$\bf\Lambda_{\bf 2}=r $ is a radius-vector of point.
We can clear up, when condition (4.1.11) is true.
The computation shows that it takes place only for
$\phi(\xi)=C\xi,
\,\beta_1=2,\,\beta_2=0,$ that is if ${\bf\Lambda^1_{\bf 1}}=
(1-CK)x+Cy,$ ${\bf\Lambda^2_{\bf 1}}= K(1-CK)x+CKy.$
The functions
$a_i(t),\,i=1,2$ can be found from the system
$$a_1(t)=-a^2_1(t)-2a_1(t)a_2(t),\quad a_1(t)=-a^2_1(t).$$

However the obtained solution has the form
${\bf u}= A(t){\bf \Lambda},\,{\bf
\Lambda}={\bf r},$ where $A'(t)=-A^2(t),$ that is we have found
a solution of form (4.1.2), with the additional condition (A1)).  It
is the solution with linear profile of velocity described before.

\subsubsection{Nonlinear transport equation on the
two-dimensional sphere.}

For the sphere of radius $r$ with the standard  orthogonal
coordinates
$x^1=\phi\in]-\pi,\pi],\,x^2=\theta\in[0,\pi]$
the metric tensor is the following:
$g_{11}=r^2 \sin^2\theta,\,g_{22}=r^2,\,
\Gamma_{11}^{2}=-\sin\theta\cos\theta,\,
\Gamma_{12}^{1}={\rm ctg}\theta,$
the other Christoffel symbols are equal to zero.

Thus, from  equation  (4.1.5) for the ratio
$z=\frac{v}{u}$ of physical components of the velocity $V=(u,v)$
we have
$$\frac{d z}{d\phi}=\cos\theta(1+z^2),\,
\frac{d \theta}{d\phi}=\sin\theta z.$$

The equation may be solved explicitly. Namely,
$$z=\pm\sqrt{-1+\sin^2\theta\Psi\left(\phi\pm \int
\limits_{\theta_0}^{\theta}\frac{d\tau}{\sin
\tau\sqrt{K+\ln\sin^2\tau}}\right)},\eqno(4.1.14)$$
with an arbitrary differential function $\Psi.$
Remark that this solution exists only in the strip
$|\sin \theta|>\exp(-\frac{K}{2}),$ it implies
$K>0.$

\medskip
{\sc Example.}
Consider at greater length the simplest case $\Psi=C=const.$
We can see from (4.1.14) that the solution is defined by the
formula only if $|\sin
\theta |\ge \frac{1}{\sqrt{C}},$
that is in some spherical strip
$\Pi$ surrounding the equator, moreover, the solution is zero
on the boundary of the strip at
$\theta_*=\pm \arcsin
\frac{1}{\sqrt{C}}.$

Further, to define the meridional component $v$
we use the equation (4.1.6), so that
$$
\frac{\partial v} {\partial \phi}+\sin\theta z\frac{\partial v}
{\partial \theta}-\cos\theta\frac{v}{z} =r\sin\theta z.$$

Further considerations we carry out for the nonnegative branch $z=$
$\sqrt{C\sin^2\theta-1}.$ For the branch with the opposite sign we
can act in a similar way. Therefore,
$$\frac{\partial v} {\partial \phi}+\sin\theta
\sqrt{C\sin^2\theta-1}\frac{\partial v} {\partial
\theta}-\frac{\cos\theta}{\sqrt{C\sin^2\theta-1}} v =r\sin\theta
\sqrt{C\sin^2\theta-1}.$$

Solving the equation we obtain
$$v=-\frac{r}{\sqrt{C}}\arcsin\frac{\sqrt{C}\cos\theta}{\sqrt{C-1}}
\frac{\sqrt{C\sin^2\theta-1}}{\sin\theta}+
\frac{\sqrt{C\sin^2\theta-1}}{\sin\theta}\Psi_1(\phi-{\mathcal
R}(\theta)),\eqno(4.1.15)$$ where $\Psi_1 $ is an arbitrary
differentiable function, and $${\mathcal R}=-\ln{\rm
tg}\frac{\theta}{2}+\frac{\sqrt{C}}{\sqrt{C-1}} {\rm
arcth}\frac{{\rm tg}^2\frac{\theta}{2}-2C+1}{2\sqrt{C}\sqrt{C-1}}.
$$

Further, for the parallel directed component we have the formula
$$u=
-\frac{r}{\sqrt{C}}\arcsin\frac{\sqrt{C}\cos\theta}{\sqrt{C-1}}
\frac{1}{\sin\theta}+ \frac{1}{\sin\theta}\Psi_1(\phi-{\mathcal
R}(\theta)).\eqno(4.1.16)$$

We see that the meridional component
$v$ vanishes on the boundary of the spherical strip
$\Pi,$ and $u$ remains bounded but, generally speaking,
does not vanish.  Therefore, the motion
takes place inside $\Pi$ independently from zones near poles of
the sphere. So the stream flows round the polar zones.

Remark that the vector field from this example is not potential.

\section{Constructing solution to the Euler
equations}

The procedure of constructing solutions to system (E1--E2)
consists of the following steps:
\vskip1cm

1. Find the velocity vector
$${\bf V}= a(t){\bf \Lambda}({\bf x}),\eqno(5.1)$$
with the $C^1$ -
smooth field $\bf\Lambda$ given on $\Sigma$
satisfying (A2).
\vskip1cm

2.
Supposing that we know $\bf V,$ find
$\rho$ and $p$ from the equations (E1), (E3), linear with respect to
these functions.
\vskip1cm

3. Choose the initial data $\rho_0$ and $p_0$ such that the
compatibility condition
$${\bf\nabla} p(t,{\bf x})=\phi(t)\rho(t,{\bf
x}){\bf\Lambda}-\rho(t,{\bf x}){\bf F}
(t,{\bf x}, a(t){\bf \Lambda},\rho, p)
\eqno(5.2)$$
holds with some function $\phi(t)$.

Remark that it may be that for a certain manifold there are no
initial conditions $\rho_0, p_0$ and the vector field $\bf\Lambda$
to satisfy (5.2).

However, as we see (section 2), there are
examples, where it really holds, at least
in the euclidean space.

If we additionally suppose that ${\bf F=F(V)}$ and it is
homogeneous with respect to $\bf V$ (the case of dry and aerodynamic
friction, and of the Coriolis force as well),
then the compatibility condition takes simpler form
$${\bf\nabla} p(t,{\bf x})=\phi(t)\rho(t,{\bf
x})({\bf\Lambda}-{\bf F
(\Lambda)}).
\eqno(5.3)$$

 %\pagebreak
\subsection{The velocity field. Integral functionals for the Euler
equations}
\subsubsection {The zero right hand side}
If  $\bf \Lambda$ is found (for example by the method described in
Section 2), then the further problem is to find the function $a(t).$

Suppose that in system (E1--E3) the mass ${\mathcal
M}=\displaystyle\int_\Sigma \rho d\Sigma, $ and the total energy
${\mathcal E}=E_k(t)+E_p(t)\,$ are conserved. Here the kinetic and
potential energies are defined as $\displaystyle E_k(t)
=\frac{1}{2}\displaystyle\int\limits_{\Sigma} \rho |V|^2 \,d\Sigma,
$ and $\displaystyle \,E_p(t)=\frac{1}{\gamma-1}
\int\limits_{\Sigma} \rho\,d\Sigma ,$ respectively. This
conservation laws always take place for the smooth solutions on the
compact manifold. If the manifold is not compact, then it is
sufficiently to require a rather quick vanishing of $\rho$ and $p$
at infinity with respect to the space variables, to fulfill the
convergence of integrals expressing the mass and energy and other
integrals, that we mention below, as well.

Recall that we seek for the velocity in the form (5.1).

Introduce the functionals
$$G_m(t)=\frac{1}{2}\int\limits_{\Sigma}
\rho|{\bf\Lambda}|^m d\Sigma,$$
$$Q(t)= \int\limits_\Sigma
\rho D d\Sigma,\quad D=\nabla_i\Lambda^i,\,m>0, $$
$$F(t)=\frac{1}{2}\int\limits_\Sigma \rho({\bf \nabla}|{\bf
\Lambda}|^2,{\bf V}) d\Sigma.$$
Denote $G_2(t):=G(t).$

\subsubsection{A potential velocity field}

Suppose that the field $\bf \Lambda$
is potential, that is $\nabla_i \Lambda^j=\nabla_j
\Lambda^i,\,i\ne j$, тR $\bf\Lambda=\frac{1}{2}{\bf \nabla}|{\bf
\Lambda}|^2.$

It is not difficult to see that here for the smooth
solutions to the system (E1--E3) the following relations hold
$$G'(t)=F(t),\,\eqno(5.1.1)$$
$$F(t)=2a(t)G(t),\,\eqno(5.1.2)$$
$$G_m'(t)=m a(t)G_m(t),\,\eqno(5.1.3)$$
$$F'(t)=2E_k(t)+Q(t),\,\eqno(5.1.4)
$$
$$E_p'(t)=-E_k'(t)=-(\gamma-1)a(t)Q(t),\,\eqno(5.1.5) $$
$$E_k(t)=  a^2(t)G(t),\,\eqno(5.1.6)$$
$$E_p(t)+a^2(t)G(t)={\mathcal E}=const.\eqno(5.1.7)$$

It follows that for the solutions with the velocity field of form
(5.1)
we obtain the system
$$G'(t)=2a(t)G(t),\,\eqno(5.1.8)$$
$$a'(t)=-a^2(t)+\frac{1}{2}\frac{Q(t)}{G(t)}.\eqno(5.1.9)$$
Moreover, from  (5.1.3) we obtain
$$G_m(t)=K_{ml} (G_l(t))^{\frac{m}{l}},\quad
K_{ml}=G_m(0)G_l(0)^{-\frac{m}{l}}=const.\eqno(5.1.10)$$

It will be more convenient to involve a new variable
$\tilde G(t)=\displaystyle\frac{1}{G(t)}>0$
and instead of (5.1.8--5.1.9) to consider  the
system
$$
a'(t)=-a^2(t)+
\frac{1}{2}{Q(t)}{\tilde G(t)},\eqno(5.1.11)$$
$$\tilde G' (t)=-2a(t)\tilde G(t).\eqno(5.1.12)$$

However this system is not closed. Below we consider the cases were
there is a possibility to obtain it in a closed form.

Note that the properties of the function
$a(t)$ are different from the case of transport equation
considered in Section 3.
\medskip

{\sc Remark 5.1.1.}
The functional $Q(t)$ equals to zero identically if $p\equiv 0$ (so
called pressure free gas dynamics or in the case of divergency free
field $\bf\Lambda$).  As follows from (5.1.9), if
$a(0)<0$,then the velocity $\bf V$ goes to infinity, and $G(t)$ goes
to zero in a finite time.  Taking into account the mass conservation
we conclude that on the points where $|{\bf\Lambda}|=0$ the mass
concentrates on the set of zero measure.  The phenomenon really takes
place in the case of the pressure free gas dynamics.  The second
possibility don't realize as $\bf\Lambda$ cannot have the zero
divergency.

\subsubsection{A non-potential velocity field}

If the field $\bf\Lambda$ satisfies (A2), but is not potential,
then the properties (5.1.1--5.1.5) take place, but the
expression (5.1.6) for the kinetic energy $E_k(t)$
becomes more complicated. For example, in the two-dimensional
space the following relation takes place:
$$E_k(t)=a^2(t)G(t)+2\int\limits_\Sigma \rho J\Lambda^1 \Lambda^2
 (1-D)\, d\Sigma,\eqno (5.1.13)$$
with $D=\nabla_i\Lambda^i,\,J=\nabla_2\Lambda^1-\nabla_1\Lambda^2.$
It is interesting that if the divergency $D$ is equal to 1,
then (5.1.13) coincide with (5.1.6) and we obtain once more the
system (5.1.11 -- 5.1.12). Recall that there are two
possible values of constant divergency, namely, $D=1$ and $D=2$ (see
Theorem 3.3), the example of a vector field $\bf\Lambda$ on the plane
with the divergency equal to 1 constructed in Subsection 2.1.1.

\subsubsection{The field $\bf \Lambda$ with a constant divergency
$D$}

Suppose once more that $\bf\Lambda$ is potential.
It is easy to see that if $D$ is constant, then
$$Q(t)=(\gamma-1){D}E_p(t),\eqno(5.1.14) $$
$$E_p(t)=const\cdot
G^{-\frac{(\gamma-1)D}{2}}(t).\eqno(5.1.15)$$

  Thus, the system (5.1.11--5.1.12) has the closed form
$$\tilde G'(t)=-2a(t)\tilde G(t),\quad \
a'(t)=-a^2(t)+ const\cdot \tilde G^{\frac{(\gamma-1)D}{2}+1}.
\eqno(5.1.16)$$

The energy conservation gives
$$const\cdot
G^{-\frac{(\gamma-1)D}{2}}(t)+a^2(t)G(t)=\mathcal E.$$ So for $D>0$
the function $G(t)$ cannot go to zero ($a(t)$ to the infinity,
respectively), therefore the mass cannot be concentrated in a point.

The only equilibrium point is the origin. The phase portrait
is presented on Figure 1.

%\begin {figure}
%\scalebox {1.35}{\includegraphics [1,150]{fig1.pdf}} \caption {}
%\end {figure}

\begin{figure}%[h!]
\centerline{\includegraphics[width=0.5\columnwidth]{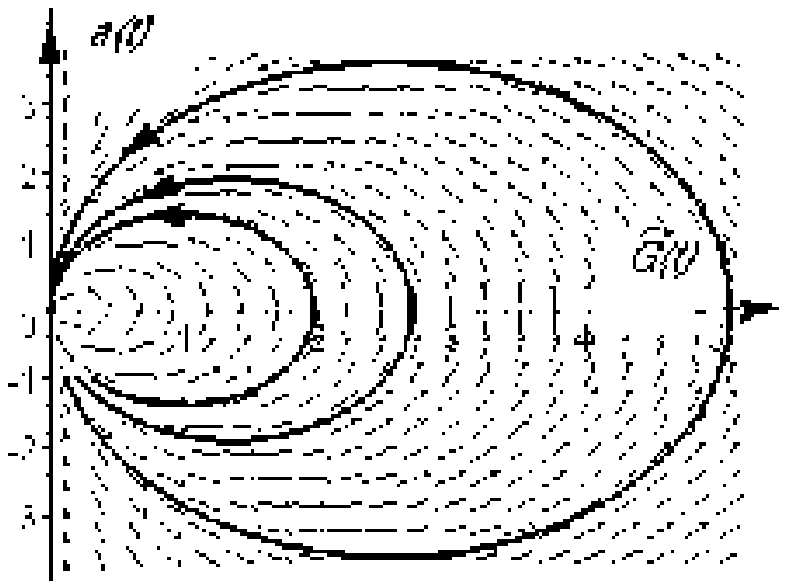}}
\caption{}
\end{figure}

\subsection{The bounded field $\bf \Lambda$ with an arbitrary
divergency $D$}

If the divergency is not constant, then, generally speaking,
we cannot obtain the closed system for finding $a(t).$
However, some properties of $a(t)$ we can get.

It is evident that $\tilde G(t)\ge 0.$ In the equilibrium points
we have $Q(t)=0.$ It may exist several points
of such kind on the axis $a=0.$ Anyway, the phase trajectories of
(5.1.11--5.1.12) are symmetric with respect to the axis.

Moreover,   the function $a(t)$ cannot go to
the infinity. Really,
from (5.1.12) we have that $\tilde G(t)=\tilde
G(0)\exp{\left(-2\int\limits_0^t a(\tau)\,d\tau\right)},$
therefore if $a(t)\to-\infty,\,t\to t_*,$ then $\tilde G(t)\to
\infty\,$ $(G(t)\to 0)\,$ $t\to t_*.$

Further we use the generalization of Lemma from \cite{Chemin}.
\newtheorem{lemma5}{Lemma 5.}
\begin{lemma5}
Let
$\Sigma$ be a smooth riemannian manifold
with  the metric tensor
$g_{ij}$ such that $|g_{ij}|$
is separated from zero by a positive
constant $g_*.$ If on the manifold there exists a function $f({\bf
x})$ such that $f({\bf x}) \ge 0,\quad f({\bf x}) \in
L_\gamma(\Sigma),\quad |{\bf x}|^2 f({\bf x})\in L_1(\Sigma),$ and
$\gamma >1,$ then
$$\int_{\Sigma} f\,d\Sigma\le C_{\gamma,
g_{ij}} \left( \int_{\Sigma} f^\gamma \,d\Sigma \right) ^
{\frac{2\gamma}{(n+2)\gamma-n}} \left( \int_{\Sigma} (|{\bf x}|^2
f\,d\Sigma \right) ^ { \frac{n(\gamma-1)}{(n+2)\gamma-n}}, $$ where
$$C_{\gamma,g_{ij}}=g_*^{\frac{n(\gamma-1)}{2((n+2)\gamma-n)}}\left(\left(
\frac{2\gamma}{n(\gamma-1)} \right) ^ {
\frac{n(\gamma-1)}{(n+2)\gamma-n}}+ \left(\frac{2\gamma}{n(\gamma-1)}
\right)^ { -\frac{2\gamma}{(n+2)\gamma-n}}\right).$$
\end{lemma5}

\vspace{1cm}

   According to (2.3.3)  the quantity $S$  remains constant along
particle  trajectories, so that
$S(t,{\bf x})\ge s_0$, where
$s_0=\displaystyle
\min\limits_{\Sigma}S_0({\bf x}).$

We use Lemma 5.1 to obtain
$$
E_p(t)
=\frac{1}{\gamma-1}\int \rho^\gamma \exp
S\,d\Sigma \ge \frac{\exp
s_0}{\gamma-1}\int \rho^\gamma \,d\Sigma
\ge \frac{K}{G^{n(\gamma-1)/2}(t)},$$
with
$$ K= \frac{1}{(\gamma-1) 2^{((n+2)\gamma-n)/2}}
\exp s_0 \left( \frac{\mathcal M}{C_{\gamma, g_{ij}}} \right) ^
{((n+2)\gamma-n)/2},$$ and $\mathcal M$ is a conserved total mass.

Therefore, if
$G(t)\to 0,\,$ $t\to t_*,$
then $E_p(t)\to \infty,$ and we obtain the contradiction with
the conservation of total energy.

\medskip

{\sc Remark 5.2.1} From the state equation (see Section 2.3)
we have $p\ge \exp{s_0}\rho^\gamma, $ where $s_0=\displaystyle
\min\limits_{\Sigma}S_0({\bf x}).$
The H$\rm\ddot o$lder inequality gives us
$${\mathcal M}^{\gamma}=\left(\int\limits_{\Sigma}\rho \,d\Sigma\right)
\le(\Omega(t))^{\gamma-1} \int\limits_{\Sigma}\rho^\gamma \,d\Sigma
\le \exp(-s_0)\sup\limits_t\Omega(t) \int\limits_{\Sigma}p\,d\Sigma
=\beta_0 E_p(t),$$
with the area of the density support  $\Omega(t)$ and the
positive constant
$\beta_0=\exp(-s_0)\sup\limits_t\Omega(t)(\gamma-1).$
Therefore,
if the
support of the  density is bounded,
the potential energy cannot vanish and there exists
an inaccessible potential energy, which cannot convert to the
kinetic one. Namely,
$$E_p(t)\ge {\mathcal M}^\gamma \beta_0^{-1}.$$
If the support of density is not bounded, then
the inaccessible energy does not exist, generally speaking
\cite{Chemin} (about the treatment on the inaccessible potential
energy see also \cite{RozZur},\cite{Roz98FAO}).

\subsubsection{Two-dimensional manifold}

This situation is more simple. Let us introduce  functionals
$$Q_m(t)=
\int\limits_{\Sigma} p D^m \,d\Sigma,\quad m=0,1,\dots,
$$
where $D=\nabla_i\Lambda^i.$
It is easy to see that $E_p(t)=(\gamma-1)Q_0(t),$
the functional $Q(t)$ introduced in Section 5.1.1
is $Q_1(t)$ in the new notation.

Further, from (E3) for the velocity of form (5.1)
we have
$$Q_m'(t)=a(t) \int\limits_{\Sigma} p (({\bf \Lambda ,\nabla}D^m)-
(\gamma-1)D^{m+1})
\,d\Sigma.\eqno(5.2.1)$$
From (3.1.5) and (3.1.10) we obtain
$({\bf \Lambda ,\nabla}D)=D-D^2+2J\,$ and  $J=D-1,$
therefore  (5.2.1) gives
$$Q_m'(t)=-m a(t) \int\limits_{\Sigma} p ((1+\frac{\gamma-1}{m})D^2
-3D+2)\,d\Sigma,\eqno(5.2.2)$$
or
$$Q_m'(t)=-m a(t) ((1+\frac{\gamma-1}{m})Q_{m+1}(t)+
3Q_m(t)-2),\quad m\in \mathbb N.\eqno(5.2.3)$$
For any $m$ the coefficient at $Q_{m+1}$ does not vanish,
therefore the system
(5.1.8), (5.1.9), (5.2.3) is also unclosed.

However it is interesting that
$${\rm sign}\, Q_m'(t)=-{\rm\, sign} \,a(t)$$ for
$m=2k+1,\,k=0,1,\dots,\,m<8(\gamma-1).$
It follows from (5.2.2), as
the expression $$(1+\frac{\gamma-1}{m})D^2 -3D+2$$
is positive and
separated from zero for  $m<8(\gamma-1).$

In particular,
$$Q_1(t)=-a(t)N(t),\eqno(5.2.4)$$ where
$\displaystyle N(t)=\int\limits_{\Sigma}p (\gamma D^2
-3D+2)p\,d\Sigma,$ moreover
$$N(t)\ge \alpha_0
E_p(t),\quad\alpha_0=\frac{(8\gamma-9)(\gamma-1)}{4\gamma},
\quad\gamma>\frac{9}{8}.$$

From (5.1.8) and (5.2.4) we have $\displaystyle\frac{dQ_1}{dG}=-
\frac{N}{2G}<0.$ Therefore there exists a decreasing function
$Q_1=Q_1(G)$ and the system
(5.1.8), (5.1.9) can be written as
$$G'(t)=2a(t)G(t),\qquad a'(t)=-a^2(t)+\frac{1}{2}\frac{Q_1(G)}{G}.
\eqno(5.2.5)$$

Let us study the equilibriums of (5.2.5). Remark that at the phase
plane $(G, a)$ the trajectories are symmetric with respect to the
axis $a=0.$ The second coordinate of the equilibriums is $a=0,$
the first one is $G$ satisfying the equation
$Q_1(G)=0.$ Recall that if the root of the equation
exists, it is unique. If the root does not exist, the first coordinate
of the equilibrium goes to the infinity, in the case if the field $
\bf \Lambda$ is unbounded. Otherwise, if $G(t)\le G_+=
const,$ it signifies that there exists no globally in time smooth
solution with the properties we try to construct.
If we use once more the variable $\tilde G={G}^{-1},$ then
we get the system
$$\tilde G'(t)=-2a(t)\tilde G(t),\qquad
a'(t)=-a^2(t)+\frac{1}{2}{\tilde Q_1(\tilde G)}{\tilde G},
$$
and the equilibrium point goes from the infinity to the origin.
The phase portrait is similar to the one from Fig.1.
If the field $\bf \Lambda$ is bounded, the restriction for $\tilde G$
has the form $\tilde G>{G_+}^{-1}.$
The situation is presented on Fig.2.

\begin{figure}%[h!]
\centerline{\includegraphics[width=0.5\columnwidth]{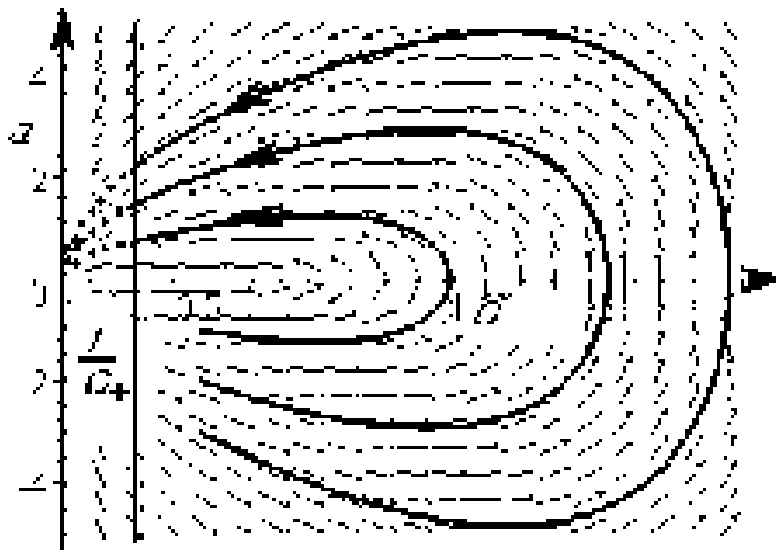}}
\caption{}
\end{figure}

%\begin {figure}
%\scalebox {1.35}{\includegraphics [1,155]{fig2.pdf}} \caption {}
%\end {figure}

If the manifold is compact, then $G(t)\le \displaystyle
\frac{1}{2}{\mathcal M}\max |{\bf \Lambda}|^2$ (if we suppose ${\bf
\Lambda}$ to be smooth), therefore there is no equilibrium in the
infinity. Note that $Q_1(t)=0$ for the constant solution
$\rho=\rho_0,\,{\bf V}={\bf 0}, \,p=p_0.$ Here the coordinates of
the only equilibrium are $(G_0,0),$ with
$G_0=\displaystyle\frac{1}{2}\int\limits_\Sigma|{\bf\Lambda}|^2\,
d\Sigma=const>0,$  it is a center (see Fig.3).

%\begin {figure}
%\scalebox {1.35}{\includegraphics [1,150]{fig3.pdf}} \caption {}
%\end {figure}

\begin{figure}%[h!]
\centerline{\includegraphics[width=0.5\columnwidth]{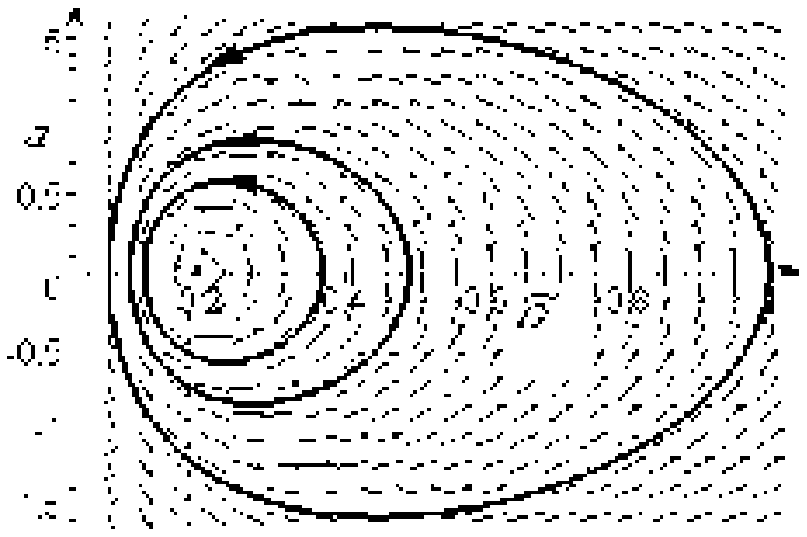}}
\caption{}
\end{figure}

If the manifold is not compact, it may exist two equilibriums -
in the finite point $(G_0,0),$ where $G_0$ is a root of the
equation $Q_1(G)=0$ (a center) and in the infinite point $(G=+\infty,
a= 0),\quad$ (the point $(\tilde G=0, a=0)$ on the plain
$(\tilde G, a)$ respectively). The situation is presented on
Fig.4.  The infinite point arises only if $\bf \Lambda$ is
unbounded.

%\begin {figure}
%\scalebox {1.35}{\includegraphics [1,150]{fig4.pdf}} \caption {}
%\end {figure}

\begin{figure}%[h!]
\centerline{\includegraphics[width=0.5\columnwidth]{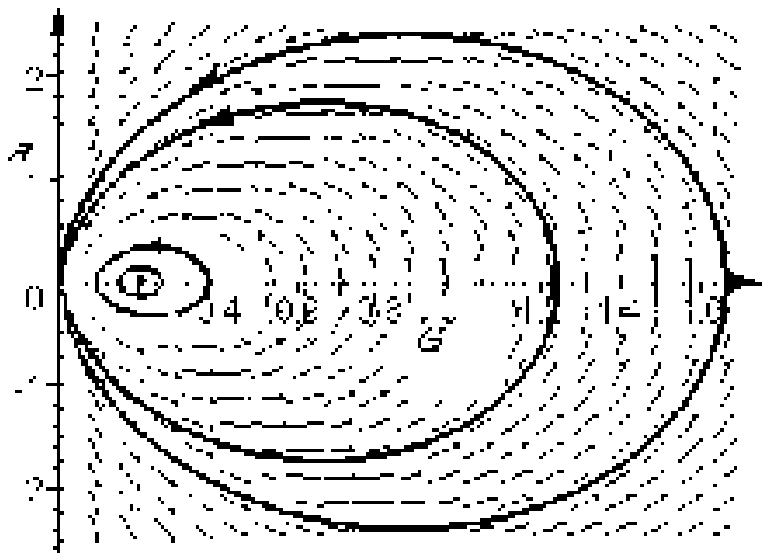}}
\caption{}
\end{figure}

{\sc Remark 5.2.2} If we suppose that there exists the limit function
$Q_\infty(t)=\displaystyle
\lim\limits_{m\to\infty}\int\limits_{\Sigma}p D^m\,d\Sigma,$
then from (5.2.3) we obtain that
the equation
$$Q_\infty(t)=-(\gamma-1)a(t)Q_\infty(t)$$
holds. Together with (5.1.8) it gives that
$$Q_\infty(t)=const\cdot G^{-\frac{\gamma-1}{2}}(t).$$
Let us stress that it is not clear whether the limit exists in the
nontrivial case $Q_\infty \ne 0.$

\subsubsection{Appearance of singularities of smooth solutions}

Suppose, for example, that  the smooth field $\bf \Lambda$ is
bounded together with its divergency $D.$ Generally speaking, the
divergency changes the sign over the manifold, that is $\inf D< 0.$
From (5.1.4) we obtain $$F'(t)=\frac{F^2(t)}{2G(t)}+Q(t)\ge
\frac{F^2(t)}{{\mathcal M}\Lambda_+^2}-(\gamma-1){\mathcal
E}D_{-},\eqno(5.2.6)$$ where
$\Lambda_+^2=\sup|\Lambda|^2,\,D_{-}=|\inf D|.$

Moreover,
$$F(t)\le 2 E_k(t)G(t)\le {\mathcal
ME}\Lambda_+^2.\eqno(5.2.7)$$

It follows from (5.2.6) that if we choose $$F(0)>
\Lambda_+\sqrt{(\gamma-1)D_{-}{\mathcal ME}},\eqno(5.2.8)$$ then
$F(t) $ will grow unboundedly, it results in a contradiction with
(5.2.7).

As the mass does not concentrate on sets of zero measure,
that $F(t)$ will grow owing to the gradient of density.
Really, in the case of a potential field
$\bf\Lambda$ we have $F(t)=\int\limits_\Sigma(\rho{\bf
V},{\bf \Lambda}) d\Sigma= -\int\limits_\Sigma{\bf \nabla} (\rho{\bf
V})\frac{|{\bf \Lambda}|^2}{2} d\Sigma=-\int\limits_\Sigma({\bf
\nabla} \rho,{\bf \Lambda}) \frac{|{\bf
\Lambda}|^2}{2}d\Sigma-\int\limits_\Sigma \rho D \frac{|{\bf
\Lambda}|^2}{2}d\Sigma.$
The last integral is bounded by the constant,
therefore, the grow of
$F(t)$ corresponds to the growth of
$|({\bf \nabla} \rho,{\bf \Lambda})|$ with respect to time.

Remark that one can ask whether the condition
(5.2.8) is possible? The comparison of (5.2.7) and (5.2.8) gives
the necessary condition for this possibility, namely,
$$D_{-}\le\frac{\mathcal ME}{\gamma-1}.$$

\subsection{The Euler equations with the right hand side of special form}

Let $\bf \Lambda$ satisfy (A2) and  the velocity field be of form
(5.1). Suppose that for ${\bf F}({\bf V})$  the following condition
holds
$$({F}({\bf V}),
{\bf\Lambda})=\sum\limits_{k=1}^{K}\lambda_k(t,a(t))|
{\bf\Lambda}|^{s_k}, \eqno(5.3.1)$$ $k\in {\mathcal N}, \,s_k\ge 0.$
Here $\lambda_k(t,a(t))$ are known functions of $a(t)$ and $t.$

\medskip
{\sc Remark 5.3.1}
The condition  (5.3.1) takes place, for example,
in the case of dry and aerodynamic friction.

Recall that for the dry friction
$${\bf F}({\bf V})=-\mu(t,{\bf x}){\bf V}.$$
Suppose $\mu(t,{\bf x})=\mu(t).$ Then
$K=1,\,s_1=2,\,\lambda_1(t,a(t))=-\mu(t)a(t).$

For the aerodynamic friction
$${\bf F}({\bf V})=-\mu_1(t,{\bf x}){\bf V}|{\bf V}|.$$
If  $\mu_1(t,{\bf x})=\mu_1(t),$
then
$K=1,\,s_1=3, \,\lambda_1(t,a(t))=-\mu_1(t)a(t)|a(t)|.$

\medskip
{\sc Remark 5.3.2}
Also (5.3.1) holds for the Navier-Stockes equation in the Euclidean
space,
where ${F}^i={\nabla}^j T^{i}_j,\,$ with the tensor $$
T^i_j=\alpha(\rho)(\nabla_i\nabla_iV^j- \nabla_j V^i)
 +\beta(\rho)\nabla_k V^k\delta^i_j,\eqno(5.3.2)$$ if
$\bf\Lambda$ coincides with the radius-vector $\bf r$ ($(\bf F({\bf
V})\equiv 0$ in the situation).

\medskip

Thus, from (5.1.1), (5.1.2) we obtain
similar to (5.1.3)
that $$F'(t)=(2a(t)G(t))'=2
E_k(t)+Q(t)+\sum\limits_{k=1}^{K}
\lambda_k(t,a(t))G_{s_k}(t),$$
and
$$a'(t)=-a^2(t)
+\frac{1}{2}\sum\limits_{k=1}^{K}
\lambda_k(t,a(t))\frac{G_{s_k}(t)}{G_2(t)}+
\frac{1}{2}\frac{Q(t)}{G_2(t)}.$$
According to (5.1.10) we have
$G_{s_k}(t)=K_{s_k}G_2^{\frac{s_k}{2}}(t)=
K_{s_k}G^{\frac{s_k}{2}}(t)$, with some
constants $K_{s_k},$ depending on the initial data.
Thus,
$$a'(t)=-a^2(t)
+\frac{1}{2}\sum\limits_{k=1}^{K}
K_{s_k 2}\lambda_k(t,a(t))G^{\frac{s_k}{2}-1}(t)+
\frac{1}{2}\frac{Q(t)}{G(t)}.$$

Introduce for the convenience (similar to Section 5.1.1) the new
variable $\tilde G(t)=G^{-1}(t),$ and obtain the system for
$(a(t),\tilde G(t), Q), \, n\in {\mathcal N}$:
$$
a'(t)=-a^2(t) +\frac{1}{2}\sum\limits_{k=1}^{K}
K_{s_k 2}\lambda_k(t,a(t))\tilde G^{1-\frac{s_k}{2}}(t)+
\frac{1}{2}{Q(t)}{\tilde G(t)},\eqno(5.3.3)$$
$$\tilde G' (t)=-2a(t)\tilde G(t).\eqno(5.3.4)$$

Generally speaking, it is unclosed, but in the case of constant
divergency $\nabla_i\Lambda^i=D,$ as in the case of zero right hand
side it is possible to express $Q(t)$ through $\tilde G(t)$ and
obtain the closed system for $(a(t), \tilde G(t)).$ Namely,
$$
a'(t)=-a^2(t) +\frac{1}{2}\sum\limits_{k=1}^{K}
K_{s_k 2}\lambda_k(t,a(t))\tilde G^{1-\frac{s_k}{2}}(t)+ \tilde K{\tilde
G}^{\frac{(\gamma-1)D+2}{2}}(t),\eqno(5.3.5)$$
$$\tilde G'
(t)=-2a(t)\tilde G(t).$$

For the  dry friction with constant friction coefficients
$\mu$  the system (5.3.5 -- 5.3.4) takes the form
$$
a'(t)=-a^2(t) -\frac{1}{2}\mu
K_{s_1 2}a(t)+ \tilde
K{\tilde G}^{\frac{(\gamma-1)D+2}{2}}(t),\eqno(5.3.6)$$
$$\tilde G'
(t)=-2a(t)\tilde G(t).$$

For the aerodynamic  friction with constant coefficients
$\mu_1$ we have
$$
a'(t)=-a^2(t) -\frac{1}{2}\mu_1 K_{s_1 2}a(t)|a(t)|\tilde
G^{-\frac{1}{2}}(t)+ \tilde K{\tilde
G}^{\frac{(\gamma-1)D+2}{2}}(t),\eqno(5.3.7)$$
$$\tilde G'
(t)=-2a(t)\tilde G(t).$$

The phase portraits of systems (5.2.6 -- 5.2.4) and (5.2.7 -- 5.2.4)
are presented on Figures
5 and 6, respectively.

%\begin {figure}
%\scalebox {1.35}{\includegraphics [1,150]{fig5.pdf}} \caption {}
%\end {figure}

%\begin {figure}
%\scalebox {1.35}{\includegraphics [1,150]{fig6.pdf}} \caption {}
%\end {figure}

\begin{figure}%[h!]
\centerline{\includegraphics[width=0.5\columnwidth]{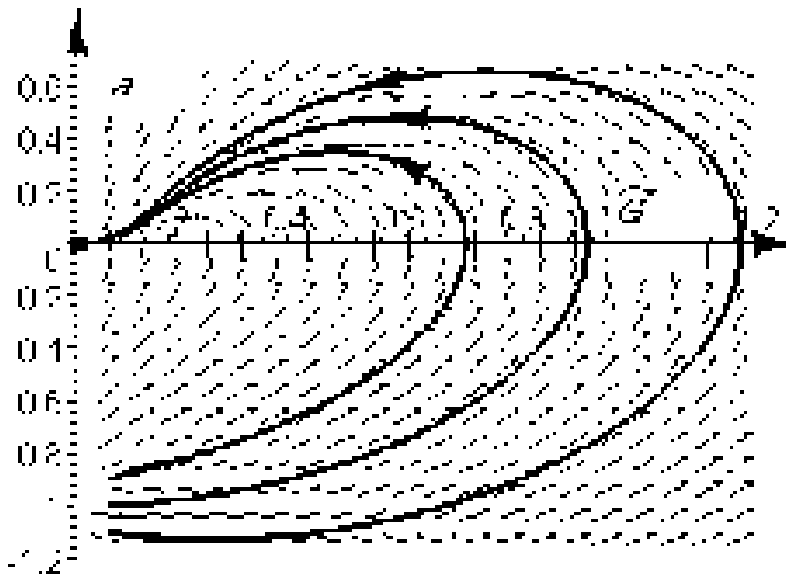}}
\caption{}
\end{figure}

\begin{figure}%[h!]
\centerline{\includegraphics[width=0.5\columnwidth]{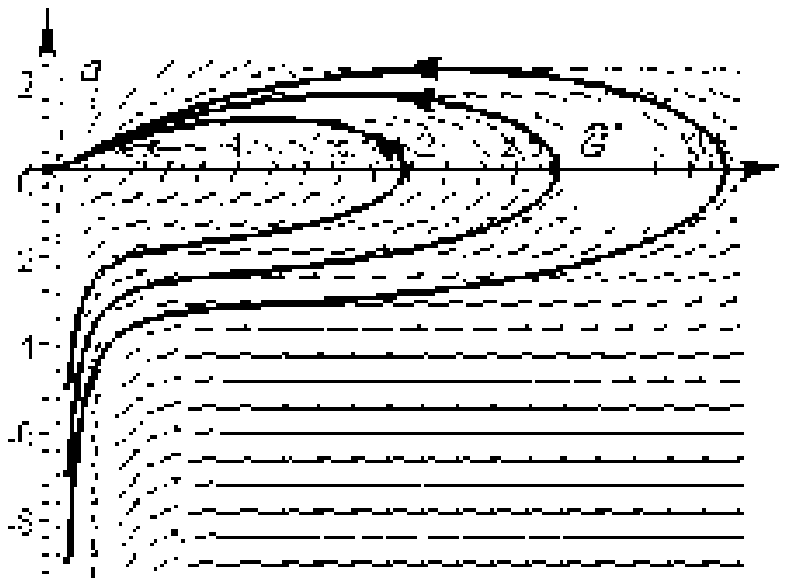}}
\caption{}
\end{figure}

In the case of dry friction there are two equilibriums:
at the origin (stable) and at the point $(0, -\frac{K_{s_1 2}\mu}{2})$
(unstable).

In the  case of aerodynamic friction (and for any friction term
with $\sigma>0$ as well (see Section 1.1)) the only equilibrium
is at the origin, moreover, it is stable.

\section{Further discussion}

There is a lot of interesting problems connecting with this
 work that would be a subject of further investigations.
Let us mention some of them.
\medskip

{\bf 1. Euler equations on the two-dimensional sphere}

a) to construct, if it is possible, a potential smooth  vector field
satisfying (A2).

b) to find a possibility to close the system (5.1.11 -- 5.1.12) and
to investigate the behaviour of the function $a(t).$

\medskip
{\bf 2. Interior solutions}

It is rather easy to show that if in the Euclidean space the vector
field $\Lambda$ satisfying (A2) differs from the radius-vector up to
the function from the Sobolev space, then the corresponding solution
is interior.

Which of solutions with velocity field (5.1) are interior?
At least is it true that any solution where $\Lambda$ has a
constant divergency is interior?

\medskip
{\bf 3. Nonviscous solutions of systems with viscous term}

It is easy to see that in the Euclidean space the solutions
with linear profile of velocity "do not feel" the viscous term
(5.3.2) in the Navier-Stokes system. Therefore the solutions
solve Navier-Stokes as well. In particular, for the viscous,
but pressure free case we can construct the solution with
the mass concentrated in the point.

It is interesting to find the solution of such kind for any
manifolds. What is a condition for the metric tensor for the
possibility of existence of "nonviscous solutions"?

Further, if a solution with linear profile of velocity is
interior for the Euler system, can we guarantee that it will
be interior for the corresponding Navier-Stokes system?

\medskip

\subsection*{Acknowledgments} The work is partially supported by
the Russian Foundation for Basic Researches Award No.
03-02-16263.

%\newpage

\end{document}